\RequirePackage{fix-cm}

\documentclass{article}

\usepackage{graphicx}
\usepackage{amssymb}
\usepackage{amsmath}
\usepackage{eucal}

\usepackage[driverfallback=dvipdfm,breaklinks=true,pdfcreator={LaTeX e Dvips}]{hyperref}
\usepackage{float}

\usepackage{mathptmx}
\usepackage[numbers]{natbib}
\usepackage{color}
\usepackage{ulem}
\usepackage{booktabs}
\usepackage{multirow}
\usepackage{textcomp}
\usepackage{bm}

\makeatletter
\let\MYcaption\@makecaption
\makeatother
\usepackage{subcaption}
      \makeatletter
\let\@makecaption\MYcaption
\makeatother

\usepackage[ruled,vlined]{algorithm2e}

\usepackage[nomarkers, nolists]{endfloat}
\renewcommand{\efloatseparator}{\vspace{1em}}

\title{Topology optimization of actively moving rigid bodies in unsteady flows}

\author{Yuta Tanabe$^\text{a,}\footnote{Corresponding author: {\tt 4524702@ed.tus.ac.jp} (Yuta Tanabe)}$,
        Kentaro Yaji$^\text{b}$ ,
        Kuniharu Ushijima$^\text{a}$ \\[12pt]
$^\text{a}$\textit{Department of Mechanical Engineering, Tokyo University of Science,}\\
              \textit{6-3-1, Niijuku,
              Katsushika-ku, Tokyo 125-8585, Japan}\\
$^\text{b}$\textit{Department of Mechanical Engineering, The University of Osaka,}\\
              \textit{2-1, Yamadaoka,
              Suita, Osaka 565-0871, Japan}}

\begin{document}

\maketitle

\begin{abstract}
    This study proposes a novel topology optimization method for unsteady fluid flows induced by actively moving rigid bodies. 
    The key idea of the proposed method is to decouple the design and analysis domains by using separate grids. 
    The design grid undergoes rigid body motion and is then overlapped onto the analysis grid. 
    After the overlap, key quantities such as the Brinkman coefficient are transferred between the grids. 
    This approach provides a direct and efficient means of representing object motion and facilitates the handling of more general and complex movements in unsteady flow conditions. 
    Since the computational cost of solving unsteady fluid problems is substantial, we employ a solver based on the lattice kinetic scheme, which is the extended version of the lattice Boltzmann method, to evaluate the design sensitivity. 
    The fundamental equations are derived, and the accuracy of the design sensitivity calculations is validated through comparison with finite difference approximations. 
    The effectiveness of the method is demonstrated through numerical examples in two-dimensional and three-dimensional settings.
    \flushleft
    \textbf{Keywords}\ \ Topology optimization $\cdot$ Fluid $\cdot$ Unsteady problem $\cdot$ Rigid body
\end{abstract}

\section{Introduction}
\label{sec1}

Everything in this world is in motion—fish and birds move gracefully through nature, while turbines and pumps operate within engineered flows.
Optimizing the shape of moving objects in fluid environments is therefore of significant interest. 
This study proposes a novel topology optimization method for unsteady fluid flows driven by moving yet non-deformable solid structures.

Topology optimization is widely regarded as one of the most effective methods in structural optimization, due to its ability to enable drastic shape variations, including the creation of internal holes.
The method was first presented in the field of structural mechanics by Bends{\o}e and Kikuchi~\cite{bendsoe1988generating}, and was later introduced in the field of fluid dynamics by Borrvall and Petersson~\cite{borrvall2003topology}.
Following their work on Stokes flow problems, subsequent studies have addressed laminar steady-state Navier-Stokes flow~\cite{gersborg2005topology}\cite{Olesen2006}, turbulent flow~\cite{kontoleontos2013adjoint}\cite{DILGEN2018363}, unsteady flow~\cite{kreissl2011topology}\cite{DENG20116688}, forced convection~\cite{yoon2010topological}\cite{Matsumori2013} and natural convection problems~\cite{alexandersen2014topology}\cite{ALEXANDERSEN2016876}.

A large number of studies on topology optimization of fluid devices have focused on stationary systems, such as flow paths~\cite{borrvall2003topology} and heat sinks~\cite{alexandersen2014topology}.
In contrast, Romero and Silva~\cite{ROMERO2014268} studied rotor optimization—an application of moving object optimization—by employing a rotating frame approach, which utilizes the Navier-Stokes equations with a centrifugal force term.
One of the key advantages of their method is that it does not require remeshing during solid motion.
As a result, the computational cost can be kept relatively low, and the method has been widely applied to three-dimensional problems~\cite{OKUBO202116}.
Subsequently, their method has been applied to turbulent flow~\cite{SA2021113551}\cite{Alonso2022}, non-Newtonian flow~\cite{Romero2017} and swirl flow~\cite{Alonso2018}\cite{ALONSO20192499}.
More recently, not only the moving parts but also the surrounding stationary components have been optimized simultaneously, such as impellers and guide vanes in diffusers~\cite{ALONSO2023592}.
The method used in the aforementioned work employs the multiple reference frame approach, in which the meshes are divided into rotating and stationary regions. 
Rotational effects are incorporated into the Navier-Stokes equations in the rotating region, and the two meshes are connected at their interface.
This method has since been applied to the design of screw-type pumps~\cite{ALONSO2025115982}.

Previous studies employing rotating coordinate systems have achieved remarkable success in the optimization of rotating machinery, partly by reducing computational cost through the use of steady-state approximations for inherently unsteady phenomena.
Although the approach is effective in the initial stage of designing devices with movable components, it is essentially difficult to apply a steady-state approximation when both movable and stationary parts are present, as in the problem settings addressed in the literature~\cite{ALONSO2023592}. 
This difficulty arises because the relative positions of the components change continuously over time—for example, between impellers and guide vanes in such configurations.
In contrast to such approaches, our proposed method directly tackles the unsteady fluid problem induced by rigid body motion, using a computationally efficient unsteady solver without relying on steady-state approximations.
This study proposes a method inspired by the Immersed Boundary Method (IBM)~\cite{PESKIN1972252}\cite{Peskin_2002} and the Volume Penalization Method (VPM)~\cite{Angot1999}\cite{BENAMOUR2020101050}. 
The key concept of the proposed approach separates the design grid, which represents the moving object, from the analysis grid, which is used for solving governing equations. 
By overlapping the design grid undergoing rigid body motion onto the analysis grid at each time step, the proposed method provides a direct representation of the object motion and enables the treatment of more general and complex movements under unsteady flow conditions.
In this study, we restrict our discussion to unsteady flows induced by rigid bodies with prescribed motion. 
However, because the representation of moving objects in the proposed framework is simple and explicit, it can be readily extended in the future to problems such as optimizing the shape of passively moving objects—like turbines or sails—as well as the control strategies themselves.

Since the fluid flow is treated as the unsteady problem in the proposed method, we have to utilize the high efficiency unsteady fluid solver.
There are several studies for unsteady fluid problems with various numerical algorithms such as the Finite Element Method (FEM)~\cite{kreissl2011topology}\cite{DENG20116688}, the Finite Volume Method (FVM)~\cite{Katsumata2025Topology}, the Spectral Element Method (SEM)~\cite{NOBIS2022105387} and the Lattice Boltzmann Method (LBM)~\cite{NORGAARD2016291}.
Although the most appropriate method depends on the specific problem to be solved, the LBM~\cite{kruger2017lattice} is a suitable and promising approach for unsteady problems in terms of computational efficiency, owing to its fully explicit algorithm and high compatibility with parallel computing.
The LBM has been applied to various three-dimensional unsteady problems, particularly in the fields of non-thermal flows~\cite{CongCHEN201717-00120}, forced convection~\cite{yaji2018}, and natural convection~\cite{Tanabe2023}.
However, the LBM requires a large amount of memory when applied to thermal flows, turbulent flows, or non-Newtonian fluids. 
In this study, the focus is placed on the representation of moving objects; therefore, the analysis is intentionally limited to laminar, isothermal, and Newtonian flow conditions, which define the scope of this paper. 
Nevertheless, with future extensibility in mind, the Lattice Kinetic Scheme (LKS)\cite{Coveney2002lattice}, a memory-efficient variant of the LBM, is employed to enable topology optimization for unsteady flow fields induced by rigid body motion.

The shape of the moving object is represented using the density method~\cite{bendsoe2003topology}. 
The design sensitivity is computed via the adjoint Lattice Kinetic Scheme (ALKS), which is an adjoint-variable method based on the LKS~\cite{tanabe2024adjointlatticekineticscheme}. 
Furthermore, to obtain binarized structures through the topology optimization procedure, we employ the Heaviside projection filter~\cite{wang2011projection}.
The shape is then updated using the Method of Moving Asymptotes (MMA)~\cite{svanberg1987method}, a gradient-based optimizer.

This paper is organized as follows:
Section~\ref{sec2} introduces the fundamental equations of the proposed method. 
The section also includes the derivation of the design sensitivity calculations by comparing them with finite difference approximations.
Section~\ref{sec3} presents numerical examples for two-dimensional and three-dimensional problems to demonstrate its effectiveness.
Finally, Section~\ref{sec4} concludes the paper by summarizing the key findings. 

\section{Formulation}
\label{sec2}

\subsection{Topology optimization for fluid}
\label{sec21}

The basic concept of the topology optimization is to regard structural optimization as the problem of defining material distribution.
In this paper, we employ the pseudo-density method, in which the material distribution is represented as the function $\gamma:D\rightarrow\left[0,1\right]$, where $0$ means fluid, $1$ means solid and interval value means porous material, respectively.
The optimization problem is expressed as:
\begin{align}
    \begin{array}{ll}
        \underset{\gamma}{\text{minimize }}\  & J(\gamma, p\left(\gamma\right), u_\alpha\left(\gamma\right))              \\
        \text{subject to }\                   & G(\gamma, p\left(\gamma\right), u_\alpha\left(\gamma\right))\leqslant 0,
        \vspace{1.5mm}                                                                                                    \\
                                              & 0\leqslant\gamma(\bm{x})\leqslant 1,\ \forall\bm{x}\in D,
    \end{array}
\end{align}
where state fields $p$ and $u_\alpha$ are governed by continuity and momentum conservation equations as follows:
\begin{align}
     & \frac{\partial u_\alpha}{\partial x_\alpha}=0, \label{eq:continuity}                                                                                                                                                      \\
     & \frac{\partial u_\alpha}{\partial t}+u_\beta\frac{\partial u_\alpha}{\partial x_\beta}=-\frac{\partial p}{\partial x_\alpha}+\nu\frac{\partial^2 u_\alpha}{\partial x_\beta^2}+F_\alpha. \label{eq:momentum_conservation}
\end{align}
Here, $\nu$ and $F_\alpha$ denote the kinematic viscosity and an external force term, respectively.
In conventional topology optimization methods for fluid flow, the external force depends on the design variable field and is formulated such that the fluid velocity approaches zero in the solid region, as proposed in the literature~\cite{borrvall2003topology}.
In our proposed method, however, this external force is modified. 
The details of this modification will be explained in Section~\ref{sec23}.
The subscripts $\alpha$ and $\beta$ denote the $x$ and $y$ directions in two-dimensional problems, and $x$, $y$ and $z$ in three-dimensional problems. 
The Einstein summation convention is applied only to repeated indices.
Equations.~\eqref{eq:continuity} and \eqref{eq:momentum_conservation} are solved in the domain $\mathcal{O}$ and the time interval $\mathcal{I}$ as an initial-boundary value problem.
The constraint functional $G$, which represents a volume constraint in this study, is given by:
\begin{equation}
    G =\frac{\int_D\gamma d\Omega}{V_{\text{max}}\int_Dd\Omega}-1, \label{eq:constraint}
\end{equation}
on the other hand, the objective functional $J$ is described in Section~\ref{sec25}.

\subsection{LKS}
\label{sec22}

\begin{figure}[t]
    \centering
    \begin{minipage}[t]{0.49\columnwidth}
        \centering
        \includegraphics[width=\columnwidth]{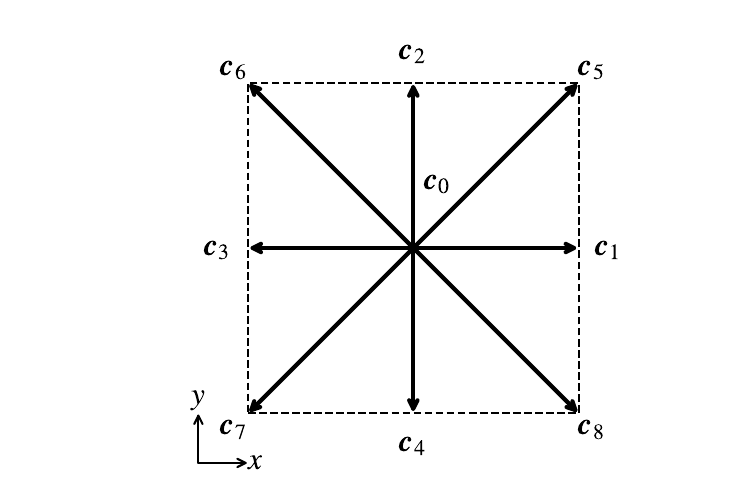}
        \subcaption{D2Q9}
        \label{fig:d2q9}
    \end{minipage}
    \begin{minipage}[t]{0.49\columnwidth}
        \centering
        \includegraphics[width=\columnwidth]{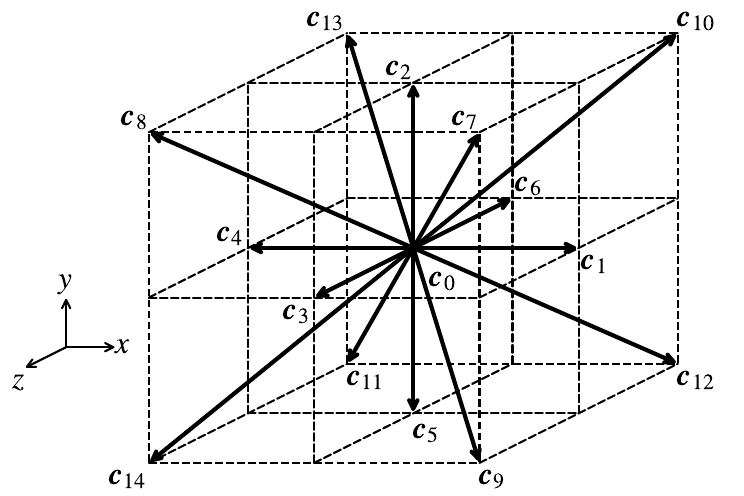}
        \subcaption{D3Q15}
        \label{fig:d3q15}
    \end{minipage}
    \caption{The velocity set for each lattice gas model.}
\end{figure}

The state fields $p$ and $u_\alpha$ are calculated by the LKS~\cite{Coveney2002lattice}, which is the special version of the LBM~\cite{kruger2017lattice}.
In the LBM, a fluid is approximated as an aggregation of fictitious particles with a finite set of velocities, and its behavior is simulated through the calculation of the motion of these particles. 
The velocity distribution functions of the fictitious particles, $f_i:\mathcal{O}\times\mathcal{I}\rightarrow\mathbb{R}$ ($i=0,1,\cdots,Q-1$), are governed by the following discrete velocity Boltzmann equations:
\begin{equation}
    \text{Sh}\frac{\partial f_i}{\partial t}+c_{i\alpha}\frac{\partial f_i}{\partial x_\alpha}=-\frac{1}{\varepsilon}\left(f_i-f_i^{\text{eq}}\right)+3\Delta xw_ic_{i\alpha}F_\alpha, \label{eq:discrete_boltzmann_equation}
\end{equation}
where $c_{i\alpha}$ and $w_i$ are the velocity set of the fictitious particles and the corresponding weights, and they are depending on the lattice gas model.
In this paper, D2Q9 model (see Fig.~\ref{fig:d2q9}) is employed for two-dimensional cases and D3Q15 model (see Fig.~\ref{fig:d3q15}) is employed for three-dimensional cases, respectively.
Here, $\text{Sh}$ denotes the Strouhal number, and $\varepsilon$ is the dimensionless parameter of the same order as the Knudsen number.
In the LBM, the lattice Boltzmann equation shown below, obtained by discretizing Eq.~\eqref{eq:discrete_boltzmann_equation}, is solved.
\begin{align}
     & f_i^*\left({\bm x}+{\bm c}_i\Delta x,t+\Delta t\right)= f_i\left({\bm x},t\right)-\frac{1}{\tau}\left\{f_i\left({\bm x},t\right)-f^\text{eq}_i\left({\bm x},t\right)\right\}, \label{eq:lattice_boltzmann_equation} \\
     & f_i\left(\bm{x},t+\Delta t\right)=f_i^*\left(\bm{x},t+\Delta t\right)+3\Delta xw_ic_{i\alpha}F_\alpha\left({\bm x},t+\Delta t\right), \label{eq:lattice_boltzmann_equation_external_force}
\end{align}
where $\tau$ is the dimensionless relaxation time related to $\varepsilon$, and is given by $\tau=\varepsilon/\Delta x$.
The fluid density $\rho$ and the fluid velocity $u_\alpha$ are obtained from the moments of the distribution function as follows:
\begin{align}
    \rho\left(\bm{x},t\right)     & =\sum_{i=0}^{Q-1}f_i\left(\bm{x},t\right),             \\
    u_\alpha\left(\bm{x},t\right) & =\sum_{i=0}^{Q-1}c_{i\alpha}f_i\left(\bm{x},t\right),
\end{align}
where the pressure is given by $p=\rho/3$. 
Hereafter, $\rho$ is used as the state variable instead of $p$.

In the LKS, the relaxation time $\tau$ in Eq.~\eqref{eq:lattice_boltzmann_equation} is set to $1$.
Consequently, the state fields without external forces are computed as:
\begin{align}
    \rho^*\left(\bm{x},t\right)     & =\sum_{i=0}^{Q-1}f_i^{\text{eq}}\left(\bm{x}-\bm{c}_i\Delta x,t-\Delta t\right), \label{eq:rho}          \\
    u_\alpha^*\left(\bm{x},t\right) & =\sum_{i=0}^{Q-1}c_{i\alpha}f_i^{\text{eq}}\left(\bm{x}-\bm{c}_i\Delta x,t-\Delta t\right), \label{eq:u}
\end{align}
where $f_i^{\text{eq}}$ is the local equilibrium distribution function, and is given by:
\begin{equation}
    f_i^{\text{eq}}=w_i\left\{\rho+3c_{i\alpha}u_\alpha+\frac{9}{2}c_{i\alpha}u_\alpha c_{i\beta}u_\beta-\frac{3}{2}u_\alpha^2+\Delta xA\left(\frac{\partial u_\alpha}{\partial x_\beta}+\frac{\partial u_\beta}{\partial x_\alpha}\right)c_{i\alpha}c_{i\beta}\right\},
\end{equation}
where $A$ is a parameter of order $O\left(1\right)$, and it is related to the kinematic viscosity as follows:
\begin{equation}
    \nu=\left(\frac{1}{6}-\frac{2}{9}A\right)\Delta x.
\end{equation}
Subsequently, the influence of external forces is incorporated through the update state fields, derived from the moments of Eq.~\eqref{eq:lattice_boltzmann_equation_external_force}:
\begin{align}
    \rho\left(\bm{x},t\right)     & =\rho^*\left(\bm{x},t\right), \label{eq:rho_ex}                                 \\
    u_\alpha\left(\bm{x},t\right) & =u_\alpha^*\left(\bm{x},t\right)+F_\alpha\left(\bm{x},t\right). \label{eq:u_ex}
\end{align}

Since the LKS only requires the macroscopic values to compute state fields, pressure or flow velocity are applied directly for the boundary conditions as follows:
\begin{align}
    \rho=\bar{\rho},u_s=\bar{u}_s & \text{ on } \Gamma_\rho, \label{eq:bc_rho} \\
    u_\alpha=\bar{u}_\alpha       & \text{ on } \Gamma_u, \label{eq:bc_u}
\end{align}
where $\Gamma_u$ and $\Gamma_\rho$ are flow velocity prescribed boundary and pressure prescribed boundary, respectively.
$\bar{\rho}$, $\bar{u}_\alpha$ and $\bar{u}_s$ are the prescribed state field values and $s$ denotes the tangential component of the flow velocity.

\subsection{Representation of moving objects}
\label{sec23}

\begin{figure}[t]
    \centering
    \includegraphics[width=\columnwidth]{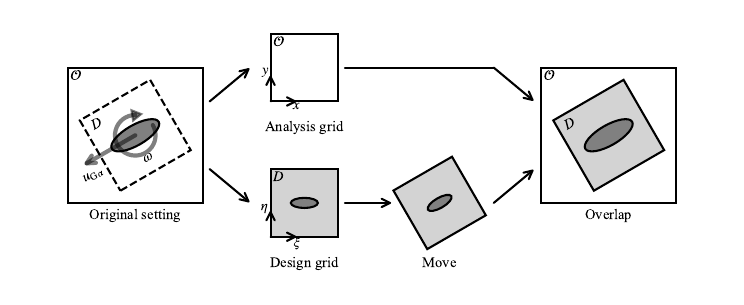}
    \caption{Representation of the moving object. The key concept is the separation between the design grid and the analysis grid. At each time step, the design grid is moved and then overlapped onto the analysis grid.
    }
    \label{fig:domains}
\end{figure}

The key concept of the proposed approach separates the grid representing the design domain $D$ from the grid representing analysis domain $\mathcal{O}$ as depicted in Fig.~\ref{fig:domains}.
While previous studies have considered limited types of rigid body motion, the proposed method accommodates arbitrary rigid body motion of the design domain.
The design domain has rigid body motion with velocity $u_{\text{G}\alpha}$ and/or angular velocity $\omega$, without any deformation.
The pseudo-density $\gamma$ is defined on $D$ and the coefficient of the Brinkman model, $\kappa^{\text{ref}}$, is calculated as follows:
\begin{equation}
    \kappa^{\text{ref}}\left(\bm{\xi}\right)=\kappa_{\text{max}}\frac{q\gamma\left(\bm{\xi}\right)}{\left\{1-\gamma\left(\bm{\xi}\right)\right\}+q}.
\end{equation}
Besides, the coordinate and moving velocity of the each point on $D$ are described by:
\begin{align}
    \bm{x}^{\text{ref}}\left(\bm{\xi},t\right)= & \left(\begin{array}{cc}
                                                                \cos{\theta\left(t\right)} & -\sin{\theta\left(t\right)} \\
                                                                \sin{\theta\left(t\right)} & \cos{\theta\left(t\right)}
                                                            \end{array}\right)\left(\begin{array}{c}\xi-\xi_{\text{G}}\\\eta-\eta_{\text{G}}\end{array}\right)+\left(\begin{array}{c}x_{\text{G}}\left(t\right)\\y_{\text{G}}\left(t\right)\end{array}\right), \label{eq:rigid_body_position}                                                         \\
    \bm{u}^{\text{ref}}\left(\bm{\xi},t\right)= & \left(\begin{array}{cc}
                                                                -\sin{\theta\left(t\right)} & -\cos{\theta\left(t\right)} \\
                                                                \cos{\theta\left(t\right)}  & -\sin{\theta\left(t\right)}
                                                            \end{array}\right)\left(\begin{array}{c}\xi-\xi_{\text{G}}\\\eta-\eta_{\text{G}}\end{array}\right)\omega\left(t\right)+\left(\begin{array}{c}u_{{\text{G}}x}\left(t\right)\\u_{{\text{G}}y}\left(t\right)\end{array}\right). \label{eq:rigid_body_velocity}
\end{align}
Here, the parameter $\kappa$, the coefficient of the Brinkman model on $\mathcal{O}$, and $u_{\text{S}\alpha}$, the moving velocity of the each point on $D$, are mapped as follows:
\begin{align}
    \kappa\left(\bm{x},t\right)             & =\int_DW\left(\bm{x},\bm{x}^{\text{ref}}\left(\bm{\xi},t\right)\right)\kappa^{\text{ref}}\left(\bm{\xi}\right)d\Omega, \label{eq:brinkman_ref}     \\
    u_{\text{S}\alpha}\left(\bm{x},t\right) & =\int_DW\left(\bm{x},\bm{x}^{\text{ref}}\left(\bm{\xi},t\right)\right)u^{\text{ref}}_\alpha\left(\bm{\xi},t\right)d\Omega, \label{eq:velocity_ref}
\end{align}
where $W$ is the weight function, also mentioned in reference~\cite{Peskin_2002} as a known function still used in the IBM, and defined as follows:
\begin{align}
    W\left(\bm{x},\bm{x}^{\text{ref}}\right) & =\prod_{\alpha=0}^{d-1}w\left(x_\alpha-x_\alpha^{\text{ref}}\right),                                                                                                    \\
    w\left(r\right)                          & =\left\{\begin{array}{lc}\frac{1}{4\Delta x}\left\{1-\cos{\left(\frac{\pi r}{2}\right)}\right\}&\left|r\right|\leq2\Delta x \\ 0&2\Delta x<\left|r\right|\end{array}\right.
\end{align}
According to the VPM, the external force is represented as follows:
\begin{equation}
    F_\alpha\left(\bm{x},t\right)=-\kappa\left(\bm{x},t\right)\left\{u_\alpha\left(\bm{x},t\right)-u_{\text{S}\alpha}\left(\bm{x},t\right)\right\}.
\end{equation}
Therefore, the flow velocity is as follows:
\begin{equation}
    u_\alpha\left(\bm{x},t\right)=\frac{u_\alpha^*\left(\bm{x},t\right)+\kappa\left(\bm{x},t\right)u_{\text{S}\alpha}\left(\bm{x},t\right)}{1+\kappa\left(\bm{x},t\right)}. \label{eq:u_ex_2}
\end{equation}

In summary, the computational procedure of the LKS with moving objects is presented in Algorithm~\ref{pc:moving_object}.
\begin{algorithm}
    \caption{LKS with moving objects procedure}
    $\rho\gets \rho_0$, $u_\alpha\gets u_{\alpha,0}$\tcp*{Initialize state fields}\
    \For{$t\gets t+\Delta t$}{
    $\bm{x}^{\text{ref}}\gets\text{Eq.~\eqref{eq:rigid_body_position}}$, $\bm{u}^{\text{ref}}\gets\text{Eq.~\eqref{eq:rigid_body_velocity}}$\tcp*{Get position and velocity of moving object}\
    $\kappa\gets\text{Eq.~\eqref{eq:brinkman_ref}}$, $u_{\text{S}\alpha}\gets\text{Eq.~\eqref{eq:velocity_ref}}$\tcp*{Get Brinkman coefficient and velocity of moving object}\
    $\rho^*\gets \text{Eq.~\eqref{eq:rho}}$, $u_\alpha^*\gets\text{Eq.~\eqref{eq:u}}$\tcp*{Get temporal values}\
    $\rho\gets\text{Eq.~\eqref{eq:rho_ex}}$, $u_\alpha\gets\text{Eq.~\eqref{eq:u_ex_2}}$\tcp*{Apply external forces and sources}\
    \If{is on boundary}{
        $\rho, u_s\gets\text{Eq.~\eqref{eq:bc_rho}}$, $u_\alpha\gets\text{Eq.~\eqref{eq:bc_u}}$\tcp*{Apply boundary conditions}\
    }
    }
    \label{pc:moving_object}
\end{algorithm}

\subsection{Sensitivity analysis}
\label{sec24}

We describe the sensitivity analysis method in the proposed approach, based on the adjoint variable method incorporated within the LKS, referred to as ALKS~\cite{tanabe2024adjointlatticekineticscheme}.
We consider a general functional $J$, which consists of the volume integral term and the boundary integral term, expressed as:
\begin{equation}
    J=\int_{\mathcal{I}}\int_\mathcal{O}J_\mathcal{O}d\Omega dt+\int_{\mathcal{I}}\int_{\partial\mathcal{O}}J_{\partial\mathcal{O}}d\Gamma dt.
\end{equation}
The governing equation is the discrete velocity Boltzmann equation.
The corresponding Lagrangian is given by:
\begin{equation}
    L=J+\int_\mathcal{I}\int_\mathcal{O}\sum_{i=0}^{Q-1}\tilde{f}_i\left\{\text{Sh}\frac{\partial f_i}{\partial t}+c_{i\alpha}\frac{\partial f_i}{\partial x_\alpha}+\frac{1}{\varepsilon}\left(f_i-f_i^{\text{eq}}\right)+3\Delta xw_ic_{i\alpha}\kappa\left(u_\alpha-u_{\text{S}\alpha}\right)\right\}d\Omega dt,
\end{equation}
where $\tilde{f}_i$ denotes the adjoint variable.
Note that, in a strict sense, the residuals of the initial and boundary conditions should also be included. 
However, we ignored these residuals to avoid notational complexity.
The derivative of the Lagrangian with respect to the design variable field $\gamma$ is given by:
\begin{align}
    \langle L^\prime,\delta\gamma\rangle= & \int_\mathcal{O}\sum_{i=0}^{Q-1}\text{Sh}\left[\tilde{f}_i\delta f_i\right]_{t_0}^{t_1}d\Omega+\int_\mathcal{I}\int_{\partial\mathcal{O}}\sum_{i=0}^{Q-1}\left(n_\alpha\left\{\tilde{f}_i\delta_{\alpha\beta}-\frac{1}{\varepsilon}\Delta xA\left(\tilde{s}_{\alpha\beta}+\tilde{s}_{\beta\alpha}\right)\right\}c_{i\beta}+\frac{\partial J_{\partial\mathcal{O}}}{\partial f_i}\right)\delta f_id\Gamma dt \notag \\                                      
                                          & +\int_\mathcal{I}\int_\mathcal{O}\sum_{i=0}^{Q-1}\left\{-\text{Sh}\frac{\partial\tilde{f}_i}{\partial t}-c_{i\alpha}\frac{\partial\tilde{f}_i}{\partial x_\alpha}+\frac{1}{\varepsilon}\left(\tilde{f}_i-\tilde{f}_i^\text{eq}\right)+3\Delta x\kappa c_{i\alpha}\tilde{u}_\alpha+\frac{\partial J_\mathcal{O}}{\partial f_i}\right\}\delta f_id\Omega dt \notag                                                   \\
                                          & +\int_\mathcal{I}\int_\mathcal{O}3\Delta x\frac{\partial\kappa}{\partial\gamma}\tilde{u}_\alpha\left(u_\alpha-u_{\text{S}\alpha}\right)\delta\gamma d\Omega dt. \label{eq:d_lagrangian}
\end{align}
Here, $\delta f_i=\left(\partial f_i/\partial\gamma\right)\delta\gamma$.
The first, second and third terms on the right-hand side in Eq.~\eqref{eq:d_lagrangian} correspond to the adjoint equation, while the last term represents the design sensitivity.

The adjoint equation is discretized and solved in the same manner as the LKS.
First, the temporal adjoint fields without any source terms $\tilde{\rho}^*$, $\tilde{u}_\alpha^*$ and $\tilde{s}_{\alpha\beta}^*$ are computed as follows:
\begin{align}
    \tilde{\rho}^*\left(\bm{x},t\right)=            & \sum_{i=0}^{Q-1}w_i\tilde{f}_i^\text{eq}\left(\bm{x}+\bm{c}_i\Delta x,t+\Delta t\right), \label{eq:arho}                    \\
    \tilde{u}_\alpha^*\left(\bm{x},t\right)=        & \sum_{i=0}^{Q-1}w_ic_{i\alpha}\tilde{f}_i^\text{eq}\left(\bm{x}+\bm{c}_i\Delta x,t+\Delta t\right), \label{eq:au}           \\
    \tilde{s}_{\alpha\beta}^*\left(\bm{x},t\right)= & \sum_{i=0}^{Q-1}w_ic_{i\alpha}c_{i\beta}\tilde{f}_i^\text{eq}\left(\bm{x}+\bm{c}_i\Delta x,t+\Delta t\right), \label{eq:as}
\end{align}
where $\tilde{f}_i^{\text{eq}}$ is corresponding to the local equilibrium function for the state fields, and described by:
\begin{equation}
    \tilde{f}_i^\text{eq}=\tilde{\rho}+3c_{i\alpha}\left(\tilde{u}_\alpha+3\tilde{s}_{\alpha\beta}u_\beta-\tilde{\rho}u_\alpha\right)-\Delta xA\frac{\partial}{\partial x_\beta}\left(\tilde{s}_{\alpha\beta}+\tilde{s}_{\beta\alpha}\right)c_{i\alpha}.
\end{equation}
Next, we incorporate source terms by solving the following equations:
\begin{align}
    \tilde{\rho}\left(\bm{x},t\right)=            & \tilde{\rho}^*\left(\bm{x},t\right)+\sum_{i=0}^{Q-1}w_i\left\{-3\Delta x\kappa\left(\bm{x},t\right)c_{i\alpha}\tilde{u}_\alpha\left(\bm{x},t\right)-\Delta x\frac{\partial J_\mathcal{O}}{\partial f_i}\left(\bm{x},t\right)\right\}, \label{eq:arho_ex}                               \\
    \tilde{u}_\alpha\left(\bm{x},t\right)=        & \tilde{u}_\alpha^*\left(\bm{x},t\right)+\sum_{i=0}^{Q-1}w_ic_{i\alpha}\left\{-3\Delta x\kappa\left(\bm{x},t\right)c_{i\beta}\tilde{u}_\beta\left(\bm{x},t\right)-\Delta x\frac{\partial J_\mathcal{O}}{\partial f_i}\left(\bm{x},t\right)\right\}, \label{eq:au_ex}                    \\
    \tilde{s}_{\alpha\beta}\left(\bm{x},t\right)= & \tilde{s}_{\alpha\beta}^*\left(\bm{x},t\right)+\sum_{i=0}^{Q-1}w_ic_{i\alpha}c_{i\beta}\left\{-3\Delta x\kappa\left(\bm{x},t\right)c_{i\gamma}\tilde{u}_\gamma\left(\bm{x},t\right)-\Delta x\frac{\partial J_\mathcal{O}}{\partial f_i}\left(\bm{x},t\right)\right\}. \label{eq:as_ex}
\end{align}
The ALKS also only requires the macroscopic values for sensitivity analysis, therefore on the boundary the adjoint variables are prescribed directly.

The functional derivative of $J$ is obtained by substituting Eq.~\eqref{eq:brinkman_ref} into the last term of Eq.~\eqref{eq:d_lagrangian}, resulting in:
\begin{align}
    \langle J^\prime,\delta\gamma\rangle & =\int_{\mathcal{I}}\int_\mathcal{O}3\Delta x\frac{\partial\kappa}{\partial\gamma}\left(u_\alpha-u_{\text{S}\alpha}\right)\tilde{u}_\alpha\delta\gamma d\Omega dt \notag                                      \\
                                         & =\int_{\mathcal{I}}\int_\mathcal{O}3\Delta x\int_DW\frac{\partial\kappa^{\text{ref}}}{\partial\gamma}d\Omega\left(u_\alpha-u_{\text{S}\alpha}\right)\tilde{u}_\alpha\delta\gamma d\Omega dt \notag           \\
                                         & =\int_{\mathcal{I}}\int_D3\Delta x\frac{\partial\kappa^{\text{ref}}}{\partial\gamma}\left(u_\alpha^{\text{ref}}-u_{\text{S}\alpha}^{\text{ref}}\right)\tilde{u}_\alpha^{\text{ref}}\delta\gamma d\Omega dt.
\end{align}

In summary, the computational procedure of the ALKS with moving objects is presented in Algorithm~\ref{pc:ALKS}.
\begin{algorithm}
    \caption{ALKS with moving objects procedure}
    $\tilde{\rho}\gets\tilde{\rho}_0$, $\tilde{u}_\alpha\gets\tilde{u}_{\alpha,0}$, $\tilde{s}_{\alpha\beta}\gets\tilde{s}_{\alpha\beta,0}$\tcp*{Initialize adjoint fields}\
    \For{$t\gets t-\Delta t$}{
    $\bm{x}^{\text{ref}}\gets\text{Eq.~\eqref{eq:rigid_body_position}}$, $\bm{u}^{\text{ref}}\gets\text{Eq.~\eqref{eq:rigid_body_velocity}}$\tcp*{Get position and velocity of moving object}\
    $\kappa\gets\text{Eq.~\eqref{eq:brinkman_ref}}$, $u_{\text{S}\alpha}\gets\text{Eq.~\eqref{eq:velocity_ref}}$\tcp*{Get Brinkman coefficient and velocity of moving object}\
    $\tilde{\rho}^*\gets\text{Eq.~\eqref{eq:arho}}$, $\tilde{u}_\alpha^*\gets\text{Eq.~\eqref{eq:au}}$, $\tilde{s}_{\alpha\beta}^*\gets\text{Eq.~\eqref{eq:as}}$\tcp*{Get temporal values}\
    $\tilde{\rho}\gets\text{Eq.~\eqref{eq:arho_ex}}$, $\tilde{u}_\alpha\gets\text{Eq.~\eqref{eq:au_ex}}$, $\tilde{s}_{\alpha\beta}\gets\text{Eq.~\eqref{eq:as_ex}}$\tcp*{Apply source terms}\
    \If{is on boundary}{
        $\tilde{\rho}\gets\bar{\tilde{\rho}}$, $\tilde{u}_\alpha\gets\bar{\tilde{u}}_\alpha$, $\tilde{s}_{\alpha\beta}\gets\bar{\tilde{s}}_{\alpha\beta}$\tcp*{Apply boundary conditions}\
    }
    }
    \label{pc:ALKS}
\end{algorithm}

\subsection{Objective functionals}
\label{sec25}

In this study, we consider two alternative objective functionals, each applied to different problem settings.
The first aims to maximize the average pressure on the boundary $\partial\mathcal{O}$, and is defined as follows:
\begin{equation}
    J_1 = -\frac{\int_I\int_{\partial\mathcal{O}}pd\Gamma dt}{\int_I\int_{\partial\mathcal{O}}d\Gamma dt}. \label{eq:objective_1}
\end{equation}
The second aims to maximize the average flow velocity in the target region $D_{\text{target}}$ given by:
\begin{align}
    J_2 = -\frac{\int_I\int_{D_{\text{target}}}n_\alpha u_\alpha d\Omega dt}{\int_I\int_{D_{\text{target}}}d\Omega dt}, \label{eq:objective_2}
\end{align}
where $n_\alpha$ denotes the unit vector in the direction of the desired flow to be maximized.
The functional $J_1$ is employed for the sensitivity verification in Section~\ref{sec26} and the numerical example in Section~\ref{sec31}, while $J_2$ is used in the examples in Sections~\ref{sec32} and \ref{sec33}.
Both objective functionals are evaluated over a time interval $I\subseteq\mathcal{I}$.

\subsection{Verification of sensitivity analysis}
\label{sec26}

\begin{figure}[t]
    \centering
    \begin{minipage}[b]{0.49\columnwidth}
        \centering
        \includegraphics[width=0.7\columnwidth]{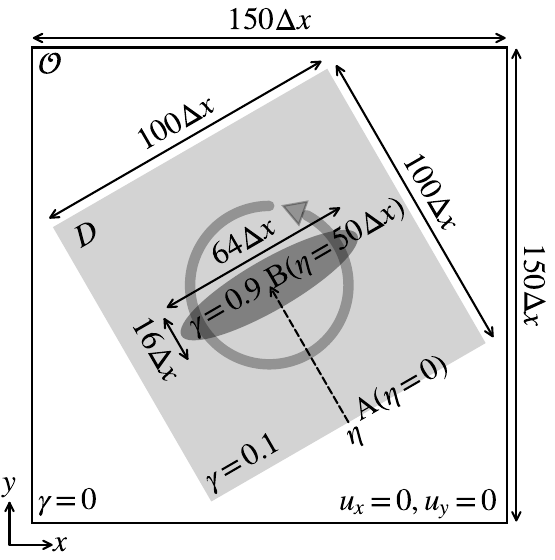}
        \subcaption{Design setting}
        \label{fig:sensitivity_design_setting}
    \end{minipage}
    \begin{minipage}[b]{0.49\columnwidth}
        \centering
        \includegraphics[width=0.9\columnwidth]{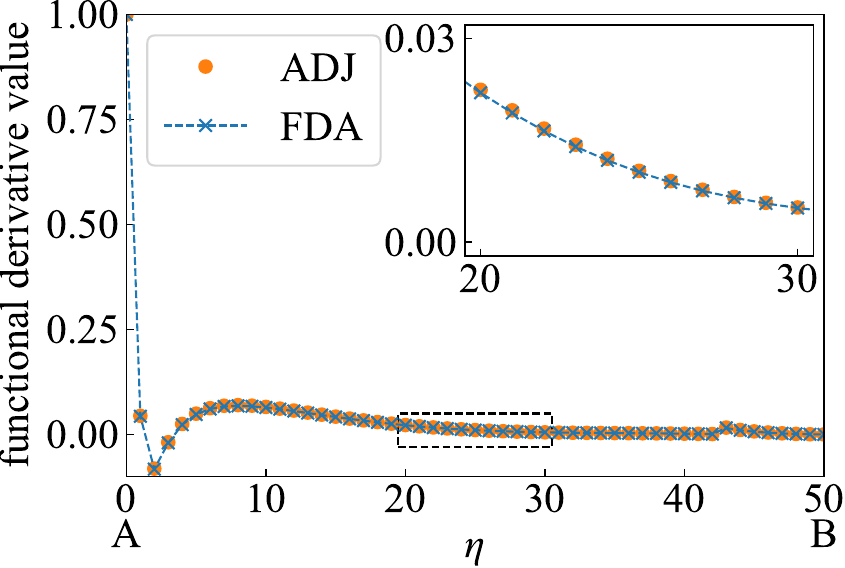}
        \subcaption{Comparison result}
        \label{fig:sensitivity_result}
    \end{minipage}
    \caption{Design setting and results of verification of sensitivity analysis. (a) The design domain rotates counterclockwise. (b) ``ADJ'' means the values by the proposed method, while ``FDA'' means ones by the FDA.}
\end{figure}

In this section, we verify the sensitivity analysis method proposed in this paper. 
The analysis domain consists of $150\Delta x\times150\Delta x$ grid points, while the design domain comprises $100\Delta x\times100\Delta x$ grid points. 
In the design domain, an elliptical region is occupied by nearly solid material $\left(\gamma=0.9\right)$, while the surrounding area is nearly fluid $\left(\gamma=0.1\right)$. 
The analysis domain is entirely fluid $\left(\gamma=0\right)$. 
The design grid rotates around the point $\left(x,y\right)=\left(75\Delta x,75\Delta x\right)$ with a time period of $3000\Delta t$. 
The values of the objective functional $J_1$, evaluated over the interval $I=\left[0,3000\Delta t\right]$ along the dashed line shown in Fig.~\ref{fig:sensitivity_design_setting}, are compared with those obtained by the Finite Difference Approximation (FDA) and shown in Fig.~\ref{fig:sensitivity_result}.
Here, ``ADJ'' means the values by the proposed method, while ``FDA'' means ones by the FDA.
Both results agree fairly with each other for any $\eta$.

\subsection{Optimization procedure}
\label{sec27}

The numerical procedure of our proposed method is written by the following steps:

\begin{enumerate}
    \item[{\it Step 1}] The design variable field $\gamma$ is initialized on the analysis domain discretized using square lattice grid.
    \item[{\it Step 2}] The state variable fields $\rho$, $u_\alpha$ are computed using the LKS, as presented in Algorithm~\ref{pc:moving_object}.
    \item[{\it Step 3}] The values of the objective functional $J$ and the constraint functional $G$ are computed.
    \item[{\it Step 4}] The adjoint variable fields $\tilde{\rho}$, $\tilde{u}_\alpha$, $\tilde{s}_{\alpha\beta}$ are computed using the ALKS, as presented in Algorithm~\ref{pc:ALKS}.
    \item[{\it Step 5}] The design sensitivities $J^\prime$ and $G^\prime$ are computed from the state and adjoint variable fields.
    \item[{\it Step 6}] The design variable field $\gamma$ is updated using the MMA~\cite{svanberg1987method}.
    \item[{\it Step 7}] The procedure returns to {\it Step 2} of the iteration loop until both the convergence criterion and the constraints are satisfied, or the maximum number of optimization steps is reached.
\end{enumerate}
The convergence criterion is defined as:
\begin{equation}
    \left|\frac{J^k-J^{k-1}}{J^{k-1}}\right|\leq 10^{-6},
\end{equation}
where the superscript $k$ denotes the optimization step. 

\section{Numerical examples}
\label{sec3}

In this section, we present numerical examples for both two-dimensional and three-dimensional moving object problems.

\subsection{2D rotor}
\label{sec31}

\begin{figure}[t]
    \centering
    \begin{minipage}[t]{0.28\columnwidth}
        \centering
        \includegraphics[width=\columnwidth]{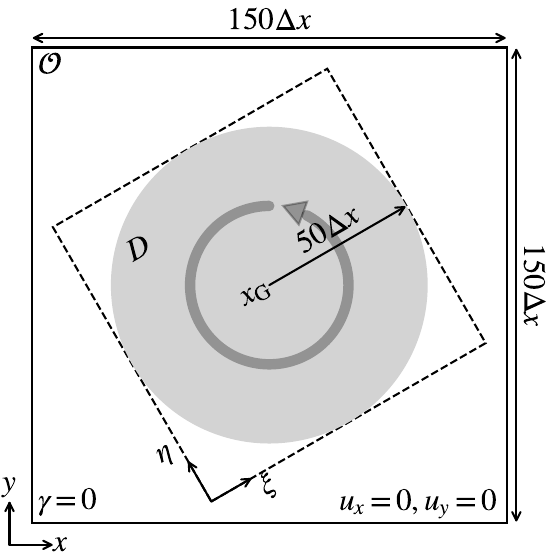}
        \subcaption{Design setting}
        \label{fig:example1_design_setting}
    \end{minipage}
    \begin{minipage}[t]{0.28\columnwidth}
        \centering
        \includegraphics[width=0.9\columnwidth]{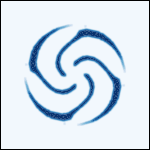}
        \subcaption{Optimized shape}
        \label{fig:example1_optimized_shape}
    \end{minipage}
    \begin{minipage}[t]{0.28\columnwidth}
        \centering
        \includegraphics[width=0.9\columnwidth]{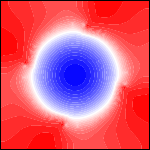}
        \subcaption{Pressure distribution}
        \label{fig:example1_pressure}
    \end{minipage}
    \begin{minipage}[t]{0.06\columnwidth}
        \centering
        \includegraphics[width=\columnwidth]{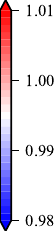}
    \end{minipage}
    \caption{Design setting, optimized shape and pressure distribution of 2D rotor. (a) The design domain rotates counterclockwise. (b) The blue area corresponds to the solid region, whereas the white area denotes the fluid region. (c) The pressure is lower in the central region and higher near the boundary.}
\end{figure}

First, we present a two-dimensional rotor design problem as an example of a rotating object in two dimensions.
The objective functional is $J_1$ defined in Eq.~\eqref{eq:objective_1}, and the constraint functional is given in Eq.~\eqref{eq:constraint}.
Here, the maximum volume fraction $V_{\text{max}}$ is set to $25\%$.
The design setting is depicted in Fig.~\ref{fig:example1_design_setting}.
The design domain $D$ rotates about the point $\left(x,y\right)=\left(75\Delta x,75\Delta x\right)$, while its center of gravity $x_{\text{G}}$ is located at $\left(\xi,\eta\right)=\left(50\Delta x,50\Delta x\right)$.
Accordingly, the position of each point on $D$ at time $t$ is described by:
\begin{equation}
    \bm{x}^{\text{ref}}\left(\bm{\xi},t\right)= \left(\begin{array}{cc}
            \cos{\left(\frac{2\pi t}{N_\text{t}}\right)} & -\sin{\left(\frac{2\pi t}{N_\text{t}}\right)} \\
            \sin{\left(\frac{2\pi t}{N_\text{t}}\right)} & \cos{\left(\frac{2\pi t}{N_\text{t}}\right)}
        \end{array}\right)\left(\begin{array}{c}\xi-50\Delta x\\\eta-50\Delta x\end{array}\right)+\left(\begin{array}{c}75\Delta x\\75\Delta x\end{array}\right),                               
\end{equation}
where $N_\text{t}$ denotes the time period, set to $3000\Delta t$.
In this example, we employ the heaviside projection filter~\cite{wang2011projection}.
The filter radius $R$ is set to $2.4\Delta x$, and the steepness parameter $\beta$ is initially set to $1$.
It is doubled every $80$ optimization steps or wherever the convergence criterion based on the fluctuation of the objective functional value is met.
\begin{figure}[t]
    \centering
    \includegraphics[width=0.5\columnwidth]{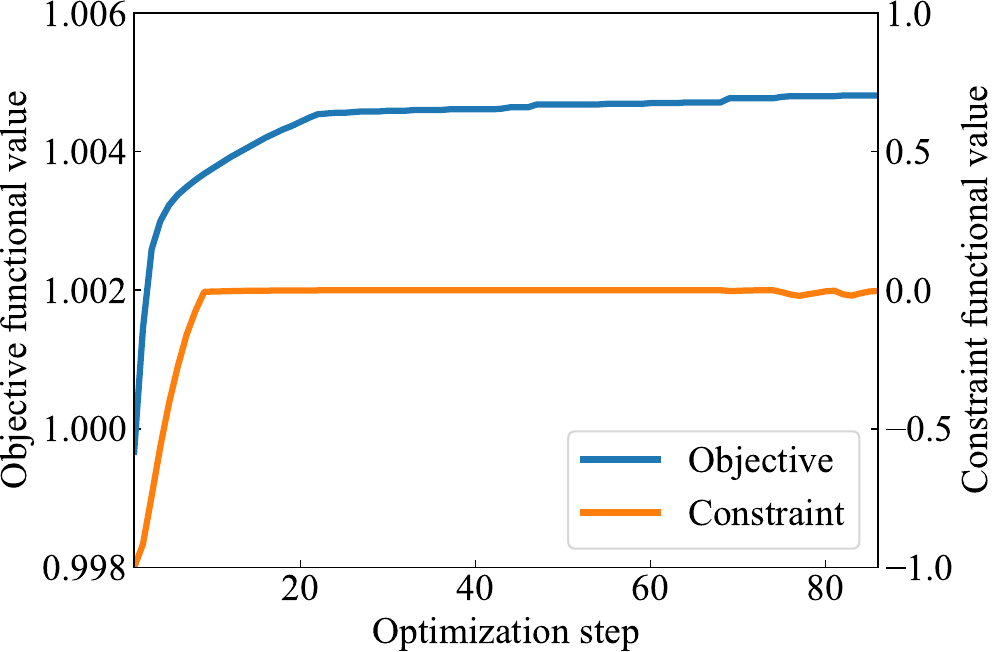}
    \caption{Convergence history of the 2D rotor. The objective functional value increases monotonically.}
    \label{fig:example1_convergence}
\end{figure}

The optimized shape is shown in Fig.~\ref{fig:example1_optimized_shape}, where the blue region represents the solid and the white region represents the fluid.
The corresponding pressure distribution is illustrated in Fig.~\ref{fig:example1_pressure}.
As shown, the pressure is lower in the central region and higher near the boundary.
The rotating solid component drives the fluid, initially located at the center, toward the outer boundary. 

The convergence history of the objective and constraint functionals is presented in Fig.~\ref{fig:example1_convergence}.
The objective functional exhibits a monotonically increasing trend throughout the optimization process.
However, abrupt changes are observed at the $47$th step and $69$th iterations, which correspond to updates of the steepness parameter $\beta$ in the projection filter.

\subsection{2D Pump}
\label{sec32}

\begin{figure}[t]
    \centering
    \begin{minipage}[t]{0.49\columnwidth}
        \centering
        \includegraphics[width=0.45\columnwidth]{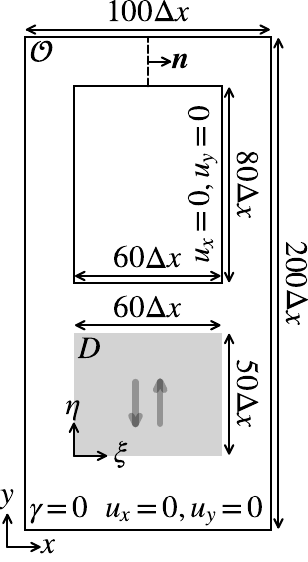}
        \subcaption{Design setting}
        \label{fig:example2_design_setting}
    \end{minipage}
    \begin{minipage}[t]{0.49\columnwidth}
\centering
        \includegraphics[width=0.4\columnwidth]{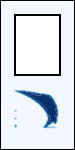}
        \subcaption{Optimized shape}
        \label{fig:example2_optimized_shape}
    \end{minipage}
    \caption{Design setting and optimized shape of translating pump. (a) The design domain undergoes periodic vertical oscillations. (b) The blue area corresponds to the solid region, whereas the white area denotes the fluid region.}
\end{figure}

\begin{figure}[t]
    \centering
    \begin{minipage}[t]{0.49\columnwidth}
        \centering
        \includegraphics[width=0.75\columnwidth]{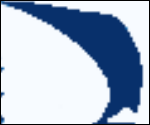}
        \subcaption{Optimized shape}
        \label{fig:example2_optimized_shape_only}
    \end{minipage}
    \begin{minipage}[t]{0.49\columnwidth}
        \centering
        \includegraphics[width=0.75\columnwidth]{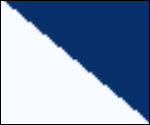}
        \subcaption{Reference shape}
        \label{fig:example2_reference_shape_only}
    \end{minipage}
    \caption{Optimized shape and reference shape of translating pump. Only the design domain is shown. The blue area corresponds to the solid region, whereas the white area denotes the fluid region.}
    \label{fig:example2_shape_only}
\end{figure}
\begin{figure}[t]
    \centering
    \begin{minipage}[t]{0.49\columnwidth}
        \centering
        \includegraphics[width=0.8\columnwidth]{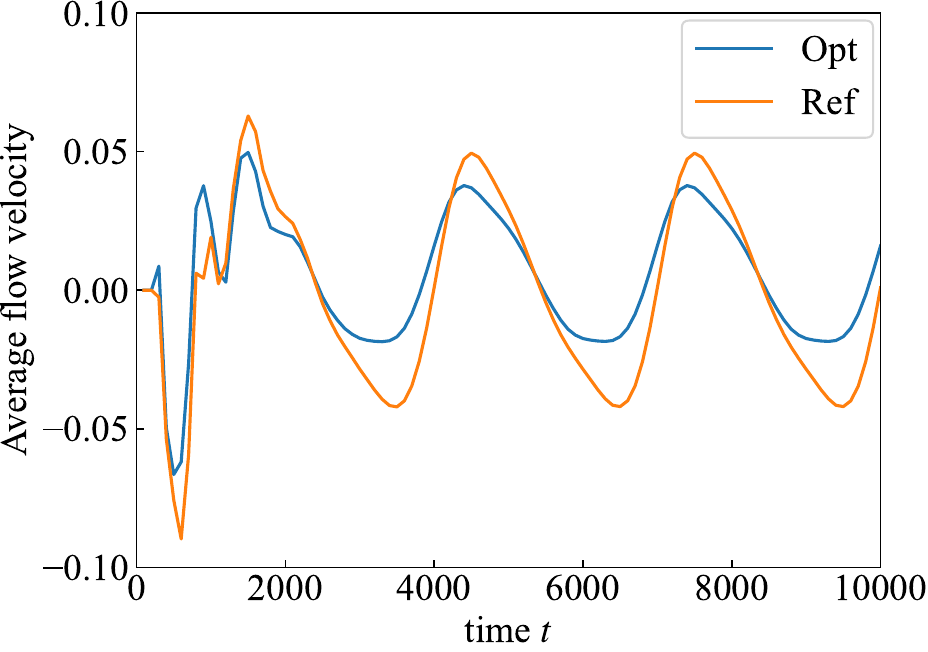}
        \subcaption{Optimized shape vs reference shape}
        \label{fig:example2_check_ref}
    \end{minipage}
    \begin{minipage}[t]{0.49\columnwidth}
        \centering
        \includegraphics[width=0.8\columnwidth]{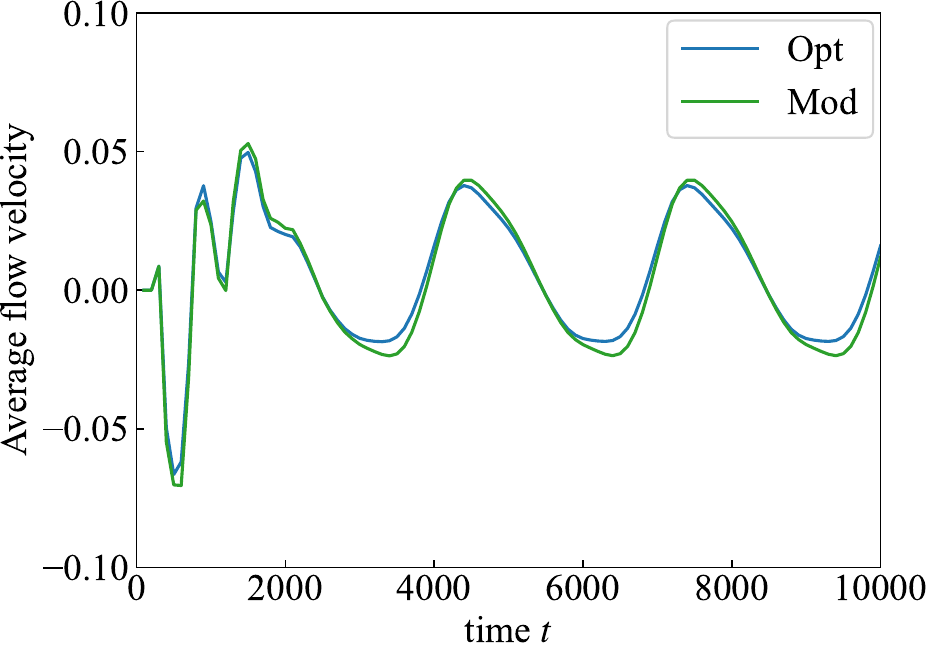}
        \subcaption{Optimized shape vs modified shape}
        \label{fig:example2_check_mod}
    \end{minipage}
    \caption{Time history of the average flow velocity on the dashed line in Fig.~\ref{fig:example2_design_setting}. The objective functional value corresponds to the time-averaged value over the quasi-steady period.}
\end{figure}

\begin{figure}[t]
    \centering
    \begin{minipage}[t]{0.24\columnwidth}
        \centering
        \includegraphics[width=\columnwidth]{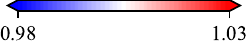}
    \end{minipage}
    \begin{minipage}[t]{0.24\columnwidth}
        \centering
        \includegraphics[width=\columnwidth]{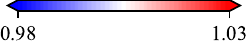}
    \end{minipage}
    \begin{minipage}[t]{0.24\columnwidth}
        \centering
        \includegraphics[width=\columnwidth]{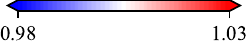}
    \end{minipage}
    \begin{minipage}[t]{0.24\columnwidth}
        \centering
        \includegraphics[width=\columnwidth]{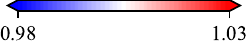}
    \end{minipage}\\
    \begin{minipage}[t]{0.24\columnwidth}
        \centering
        \includegraphics[width=\columnwidth]{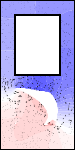}
        \subcaption{$t=0/8N_t$ (Opt)}
    \end{minipage}
    \begin{minipage}[t]{0.24\columnwidth}
        \centering
        \includegraphics[width=\columnwidth]{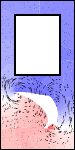}
        \subcaption{$t=1/8N_t$ (Opt)}
    \end{minipage}
    \begin{minipage}[t]{0.24\columnwidth}
        \centering
        \includegraphics[width=\columnwidth]{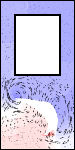}
        \subcaption{$t=2/8N_t$ (Opt)}
    \end{minipage}
    \begin{minipage}[t]{0.24\columnwidth}
        \centering
        \includegraphics[width=\columnwidth]{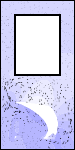}
        \subcaption{$t=3/8N_t$ (Opt)}
    \end{minipage}\\
    \begin{minipage}[t]{0.24\columnwidth}
        \centering
        \includegraphics[width=\columnwidth]{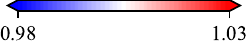}
    \end{minipage}
    \begin{minipage}[t]{0.24\columnwidth}
        \centering
        \includegraphics[width=\columnwidth]{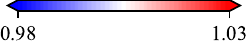}
    \end{minipage}
    \begin{minipage}[t]{0.24\columnwidth}
        \centering
        \includegraphics[width=\columnwidth]{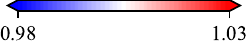}
    \end{minipage}
    \begin{minipage}[t]{0.24\columnwidth}
        \centering
        \includegraphics[width=\columnwidth]{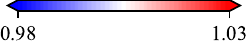}
    \end{minipage} \\
    \begin{minipage}[t]{0.24\columnwidth}
        \centering
        \includegraphics[width=\columnwidth]{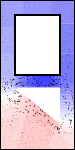}
        \subcaption{$t=0/8N_t$ (Ref)}
    \end{minipage}
    \begin{minipage}[t]{0.24\columnwidth}
        \centering
        \includegraphics[width=\columnwidth]{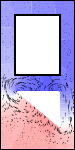}
        \subcaption{$t=1/8N_t$ (Ref)}
    \end{minipage}
    \begin{minipage}[t]{0.24\columnwidth}
        \centering
        \includegraphics[width=\columnwidth]{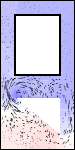}
        \subcaption{$t=2/8N_t$ (Ref)}
    \end{minipage}
    \begin{minipage}[t]{0.24\columnwidth}
        \centering
        \includegraphics[width=\columnwidth]{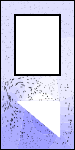}
        \subcaption{$t=3/8N_t$ (Ref)}
    \end{minipage}
    \caption{Pressure distribution and velocity vector plot for 2D pump at each time ($0/8N_t\le t\le 3/8N_t$). Here, $N_t = 3000\Delta t$.}
    \label{fig:example2_distribution_1}
\end{figure}
\begin{figure}[t]
    \centering
    \begin{minipage}[t]{0.24\columnwidth}
        \centering
        \includegraphics[width=\columnwidth]{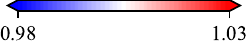}
    \end{minipage}
    \begin{minipage}[t]{0.24\columnwidth}
        \centering
        \includegraphics[width=\columnwidth]{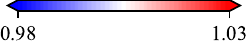}
    \end{minipage}
    \begin{minipage}[t]{0.24\columnwidth}
        \centering
        \includegraphics[width=\columnwidth]{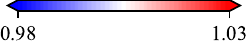}
    \end{minipage}
    \begin{minipage}[t]{0.24\columnwidth}
        \centering
        \includegraphics[width=\columnwidth]{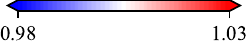}
    \end{minipage}\\
    \begin{minipage}[t]{0.24\columnwidth}
        \centering
        \includegraphics[width=\columnwidth]{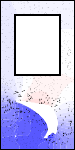}
        \subcaption{$t=4/8N_t$ (Opt)}
    \end{minipage}
    \begin{minipage}[t]{0.24\columnwidth}
        \centering
        \includegraphics[width=\columnwidth]{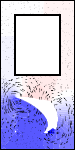}
        \subcaption{$t=5/8N_t$ (Opt)}
    \end{minipage}
    \begin{minipage}[t]{0.24\columnwidth}
        \centering
        \includegraphics[width=\columnwidth]{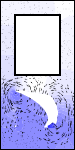}
        \subcaption{$t=6/8N_t$ (Opt)}
    \end{minipage}
    \begin{minipage}[t]{0.24\columnwidth}
        \centering
        \includegraphics[width=\columnwidth]{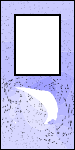}
        \subcaption{$t=7/8N_t$ (Opt)}
    \end{minipage}\\
    \begin{minipage}[t]{0.24\columnwidth}
        \centering
        \includegraphics[width=\columnwidth]{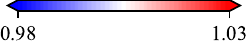}
    \end{minipage}
    \begin{minipage}[t]{0.24\columnwidth}
        \centering
        \includegraphics[width=\columnwidth]{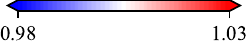}
    \end{minipage}
    \begin{minipage}[t]{0.24\columnwidth}
        \centering
        \includegraphics[width=\columnwidth]{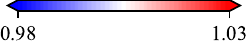}
    \end{minipage}
    \begin{minipage}[t]{0.24\columnwidth}
        \centering
        \includegraphics[width=\columnwidth]{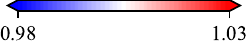}
    \end{minipage} \\
    \begin{minipage}[t]{0.24\columnwidth}
        \centering
        \includegraphics[width=\columnwidth]{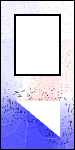}
        \subcaption{$t=4/8N_t$ (Ref)}
    \end{minipage}
    \begin{minipage}[t]{0.24\columnwidth}
        \centering
        \includegraphics[width=\columnwidth]{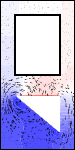}
        \subcaption{$t=5/8N_t$ (Ref)}
    \end{minipage}
    \begin{minipage}[t]{0.24\columnwidth}
        \centering
        \includegraphics[width=\columnwidth]{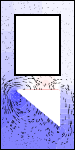}
        \subcaption{$t=6/8N_t$ (Ref)}
    \end{minipage}
    \begin{minipage}[t]{0.24\columnwidth}
        \centering
        \includegraphics[width=\columnwidth]{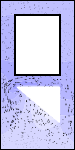}
        \subcaption{$t=7/8N_t$ (Ref)}
    \end{minipage}
    \caption{Pressure distribution and velocity vector plot for 2D pump at each time ($4/8N_t\le t\le 7/8N_t$). Here, $N_t = 3000\Delta t$.}
    \label{fig:example2_distribution_2}
\end{figure}

Next, we present a two-dimensional pump design problem as an example of translational motion in two dimensions.
The objective functional is $J_2$ defined in Eq.~\eqref{eq:objective_2}, where $D_\text{target}$ correspond to the dashed line at the top of the domain.
The constraint functional is given in Eq.~\eqref{eq:constraint}, with the maximum volume fraction $V_{\text{max}}$ set to $25\%$.
The design setting is illustrated in Fig.~\ref{fig:example2_design_setting}.
The analysis domain $\mathcal{O}$ includes an inner wall described by the region $\left\{\left(x,y\right)|20\Delta x\leq x\leq80\Delta x\cap 100\Delta x\leq y\leq180\Delta x\right\}$.
The design domain, shown in gray in Fig.~\ref{fig:example2_design_setting}, undergoes periodic vertical translation.
Accordingly, the position of each point on $D$ at time $t$ is given by:
\begin{equation}
    \bm{x}^{\text{ref}}\left(\bm{\xi},t\right)=\left(\begin{array}{cc}
            1 & 0  \\
            0 & 1
        \end{array}\right)\left(\begin{array}{c}\xi-30\Delta x\\\eta-25\Delta x\end{array}\right)+\left(\begin{array}{c}50\Delta x\\25\Delta x\sin{\left(\frac{2\pi t}{N_\text{t}}\right)+50\Delta x}\end{array}\right), \label{eq:rigid_body_position_ex2}
\end{equation}
where $N_\text{t}$ represents the time period and is set to $3000\Delta t$.
As in the previous example, we employ the heaviside projection filter, using the same parameter values.
Unlike the previous case, we instead initialize the state field values of each optimization step using the terminal values from the previous step.
The same treatment is applied to the adjoint fields, in accordance with the method proposed in a previous study~\cite{Tanabe2023} to accommodate periodic flow behavior.  
The optimized shape is shown in Fig.~\ref{fig:example2_optimized_shape}.

We compare the objective functional value of the optimized shape with that of a reference shape, which is shown in Fig.~\ref{fig:example2_reference_shape_only}.
It is worth noting that the shapes in Figs.\ref{fig:example2_shape_only} exhibit clear solid-fluid interfaces, whereas the shape in Fig.\ref{fig:example2_optimized_shape} shows a blurred interface.
This difference arises because Figs.\ref{fig:example2_shape_only} illustrate the design grid, while Fig.\ref{fig:example2_optimized_shape} illustrates the analysis grid.
The design variable is defined on the design grid and then mapped onto the analysis grid in the same manner as described in Eq.\eqref{eq:brinkman_ref}.
The average flow rate for both the optimized and reference shapes are presented in Fig.~\ref{fig:example2_check_ref}, where ``Opt'' and ``Ref'' indicate the values corresponding to the optimized and reference shapes, respectively. 
The time region $0\leq t\leq2500\Delta t$ is transient phase, while $t\leq2500\Delta t$ corresponds to the quasi-steady phase.
The value of the objective functional $J_2$ is defined as the time-averaged value over one period of the quasi-steady interval shown in Fig.~\ref{fig:example2_check_ref}.
The values obtained are $J_2=3.28\times10^{-4}$ for the optimized shape and $J_2=3.55\times10^{-5}$ for the reference shape.
The improved performance of the optimized shape can be attributed to its ability to suppress reverse flow.

Figs.~\ref{fig:example2_distribution_1} and \ref{fig:example2_distribution_2} show the pressure distribution and the fluid velocity vector plot for both the optimized and reference shapes at each time step.
The upper surface of the optimized shape is relatively smooth, resulting in a milder pressure gradient and reduced reverse flow.
In contrast, the reference shape exhibits nearly symmetric flow behavior along the flow direction, leading to less efficient performance.

In the optimized shape, three isolated solid components are present in the lower-left region of the design domain.
A comparative analysis was performed to investigate their effect.
Figure.~\ref{fig:example2_check_mod} shows the average flow rate at each time, where ``Opt'' and ``Mod'' correspond to the shapes with and without these solid island, respectively.
The objective functional value slightly decreases when the islands are removed; therefore, the islands play an important role, although their contribution is not dominant.
It is believed that the solid islands help prevent the formation of a large vortex near the lower-left region of the design domain.
Although the specific shapes differ, the reference shape—which lacks solid islands—tends to generate a large vortex in that region, whereas the optimized shape suppresses its formation.
By suppressing vortex formation, the optimized shape prevents a significant pressure drop in this region, which in turn helps to avoid undesirable backflow.

\subsection{3D rotor}
\label{sec33}

\begin{figure}[t]
    \centering
    \begin{minipage}[t]{0.32\columnwidth}
        \centering
        \includegraphics[width=\columnwidth]{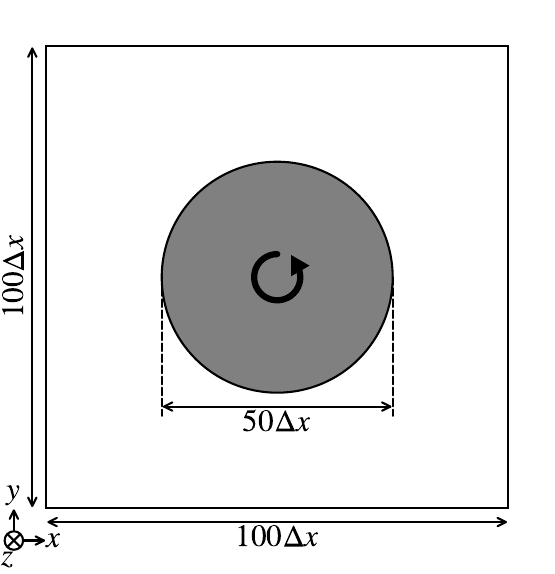}
        \subcaption{Bottom view}
    \end{minipage}
    \begin{minipage}[t]{0.32\columnwidth}
        \centering
        \includegraphics[width=\columnwidth]{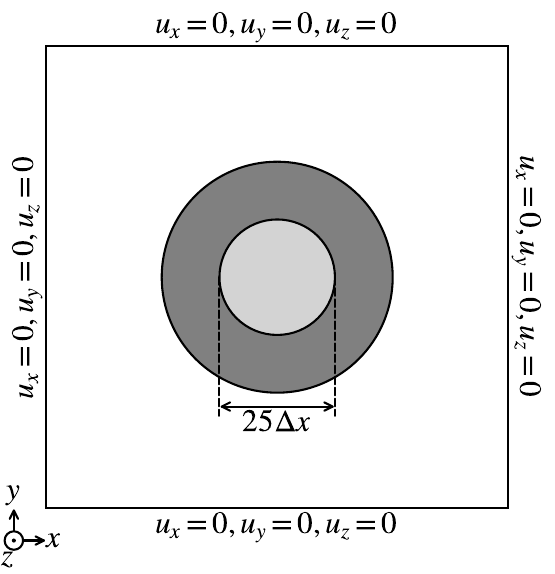}
        \subcaption{Top view}
    \end{minipage}
    \begin{minipage}[t]{0.32\columnwidth}
        \centering
        \includegraphics[width=\columnwidth]{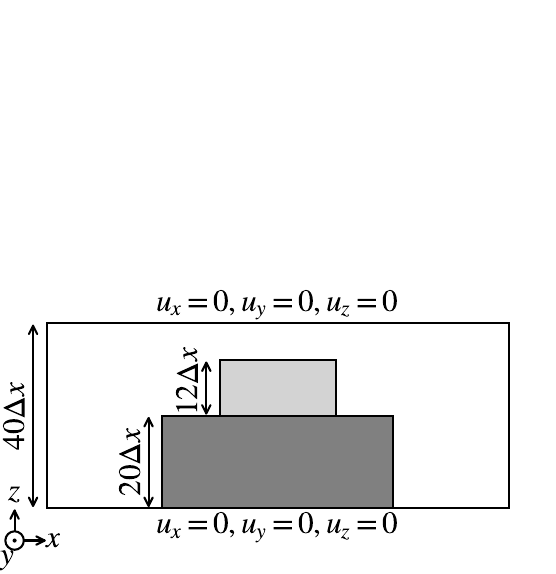}
        \subcaption{Side view}
    \end{minipage}
    \caption{Design setting of 3D rotor. The design domain rotates at the bottom of the analysis domain.}
    \label{fig:example3_design_setting}
\end{figure}

Finally, we present a three-dimensional rotor design problem to demonstrate the scalability of the proposed method.
The design setting is illustrated in Fig.~\ref{fig:example3_design_setting}.
The objective functional is $J_2$ defined in Eq.~\eqref{eq:objective_2}, where $D_\text{target}$ is a cylindrical region located at the top of the design domain, with a diameter and height of $25\Delta x$ and $12\Delta x$, respectively.
The constraint functional is given in Eq.~\eqref{eq:constraint}, with the maximum volume fraction $V_\text{max}$ set to $40\%$.
It is important to note that, in contrast to the two-dimensional cases, the volume constraint $V_\text{max}$ is initially set to $100\%$ (i.e., no constraint), and then reduced to $40\%$ in order to avoid convergence to a local optimum.
This strategy will be discussed in more detail later.
The design domain $D$ rotates about a vertical axis located at $\left(x,y\right)=\left(50\Delta x,50\Delta x\right)$.
Accordingly, the position of each point on $D$ at time $t$ is expressed as:
\begin{equation}
    \bm{x}^{\text{ref}}\left(\bm{\xi},t\right)= \left(\begin{array}{ccc}
            \cos{\left(\frac{2\pi t}{N_\text{t}}\right)} & -\sin{\left(\frac{2\pi t}{N_\text{t}}\right)} & 0  \\
            \sin{\left(\frac{2\pi t}{N_\text{t}}\right)} & \cos{\left(\frac{2\pi t}{N_\text{t}}\right)}  & 0  \\
            0                                            & 0                                             & 1
        \end{array}\right)\left(\begin{array}{c}\xi-50\Delta x\\\eta-50\Delta x\\\zeta-10\Delta x\end{array}\right)+\left(\begin{array}{c}75\Delta x\\75\Delta x\\10\Delta x\end{array}\right),                               
\end{equation}
where $N_t$ denotes the time period, which is set to $3000\Delta t$.
As in the previous two examples, the heaviside projection filter is employed, with the filter radius $R$ set to $2.4\Delta x$.
During the first $80$ optimization steps, the steepness parameter $\beta$ is set to $1$ and is subsequently doubled every $40$ steps or wherever the fluctuation of the objective functional value converges.
Similar to the previous example, both the state and adjoint field values from the terminal time of the previous optimization step are reused as the initial values for the current step.

\begin{figure}[t]
    \centering
    \begin{minipage}[t]{0.49\columnwidth}
        \includegraphics[width=\columnwidth]{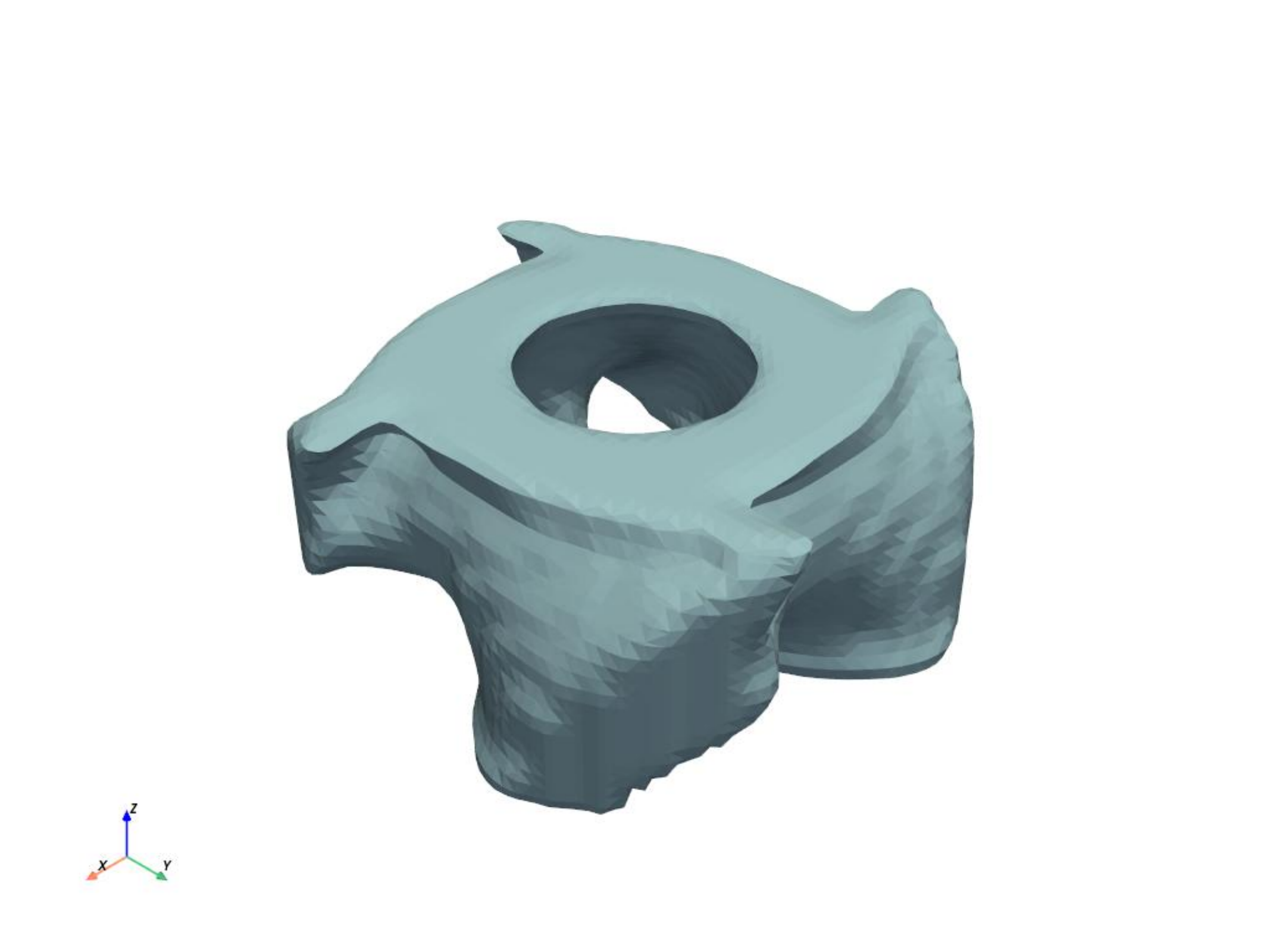}
        \subcaption{Optimized shape}
        \label{fig:example3_optimized_shape}
    \end{minipage}
    \begin{minipage}[t]{0.49\columnwidth}
        \includegraphics[width=\columnwidth]{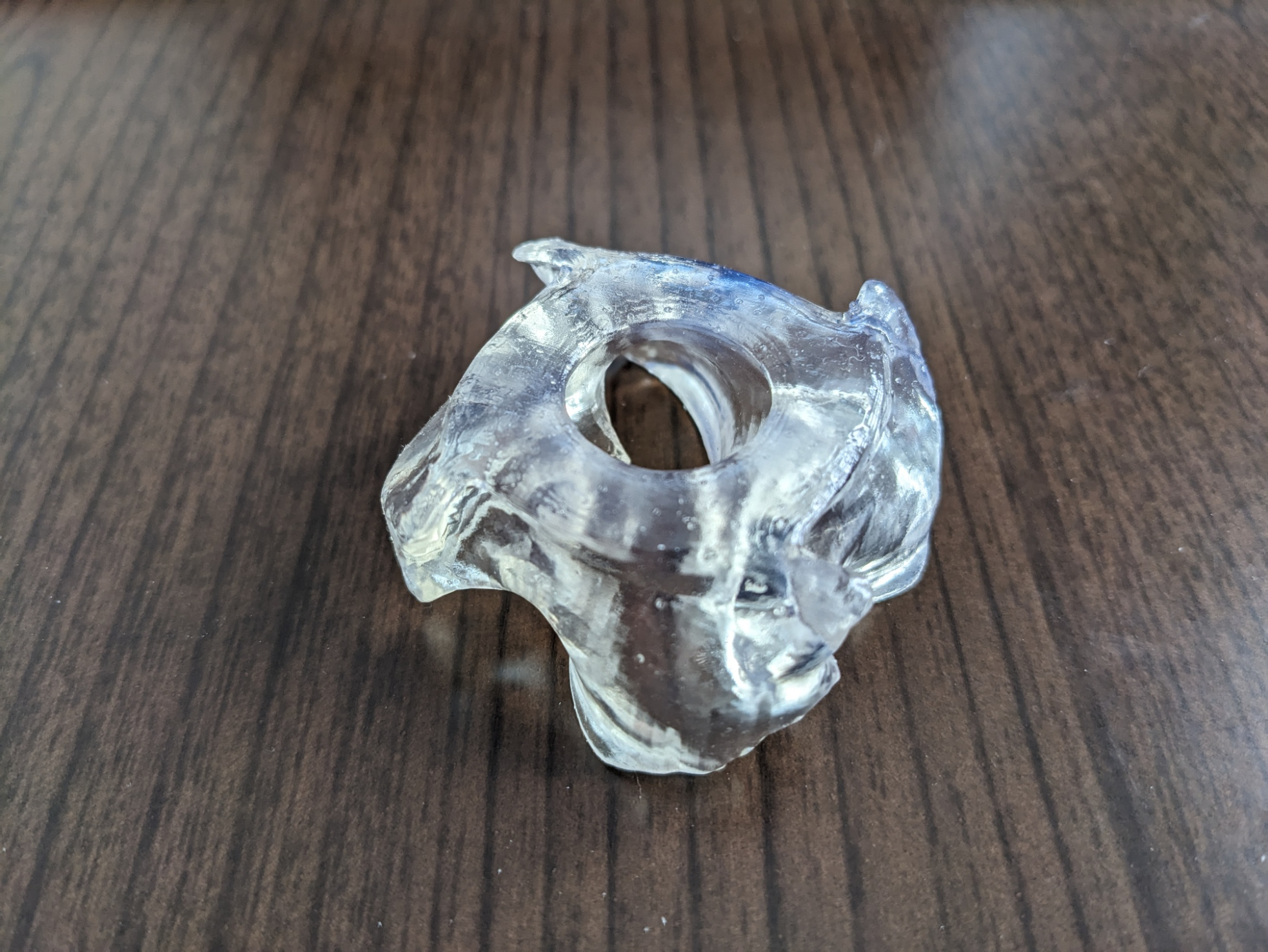}
        \subcaption{Prototype manufactured}
        \label{fig:example3_optimized_shape_printed}
    \end{minipage}
    \caption{Optimized shape and prototype manufactured of 3D rotor}
\end{figure}
\begin{figure}[t]
    \centering
    \begin{minipage}[t]{0.4\columnwidth}
        \centering
        \includegraphics[width=\columnwidth]{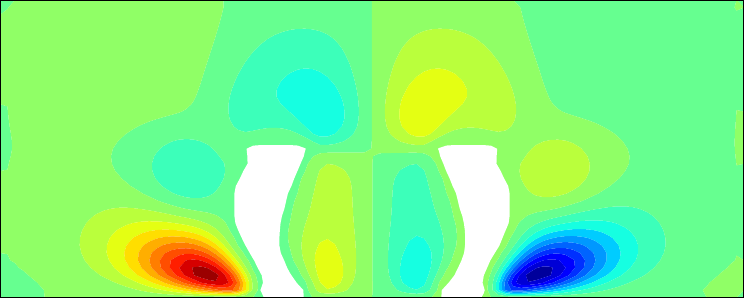}
        \subcaption{$x$ component at $t=0/3N_t$}
    \end{minipage}
    \begin{minipage}[t]{0.06\columnwidth}
        \centering
        \includegraphics[width=\columnwidth]{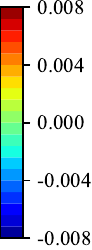}
    \end{minipage}
    \begin{minipage}[t]{0.4\columnwidth}
        \centering
        \includegraphics[width=\columnwidth]{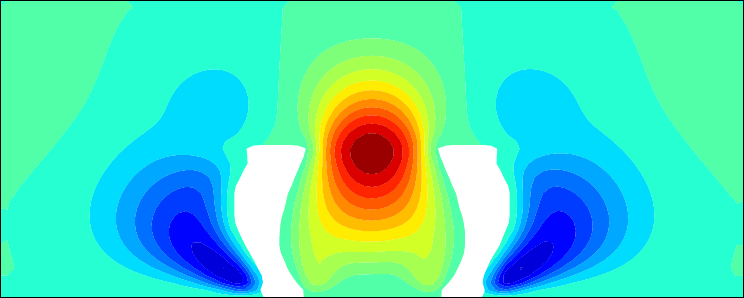}
        \subcaption{$z$ component at $t=0/3N_t$}
        \label{fig:example3_flow_x_t1000}
    \end{minipage}
    \begin{minipage}[t]{0.06\columnwidth}
        \centering
        \includegraphics[width=\columnwidth]{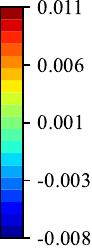}
    \end{minipage} \\
    \begin{minipage}[t]{0.4\columnwidth}
        \centering
        \includegraphics[width=\columnwidth]{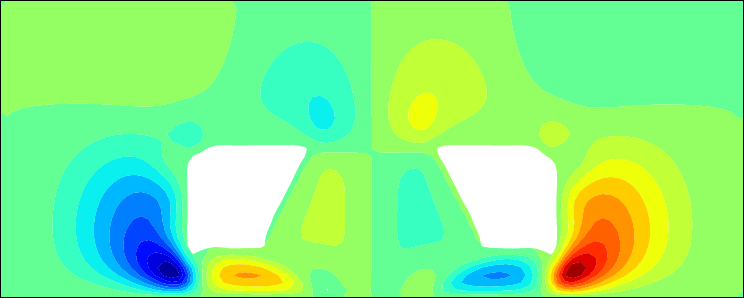}
        \subcaption{$x$ component at $t=1/3N_t$}
    \end{minipage}
    \begin{minipage}[t]{0.06\columnwidth}
        \centering
        \includegraphics[width=\columnwidth]{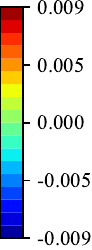}
    \end{minipage}
    \begin{minipage}[t]{0.4\columnwidth}
        \centering
        \includegraphics[width=\columnwidth]{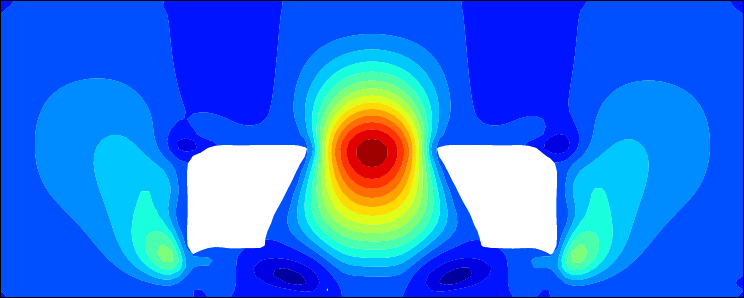}
        \subcaption{$z$ component at $t=1/3N_t$}
    \end{minipage}
    \begin{minipage}[t]{0.06\columnwidth}
        \centering
        \includegraphics[width=\columnwidth]{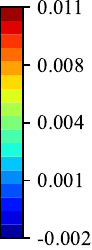}
    \end{minipage} \\
    \begin{minipage}[t]{0.4\columnwidth}
        \centering
        \includegraphics[width=\columnwidth]{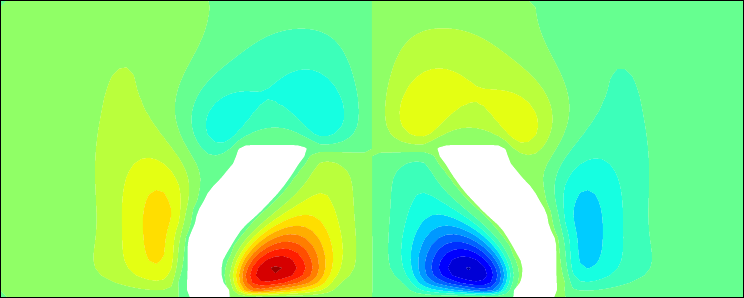}
        \subcaption{$x$ component at $t=2/3N_t$}
    \end{minipage}
    \begin{minipage}[t]{0.06\columnwidth}
        \centering
        \includegraphics[width=\columnwidth]{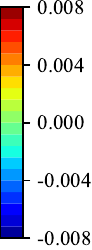}
    \end{minipage}
    \begin{minipage}[t]{0.4\columnwidth}
        \centering
        \includegraphics[width=\columnwidth]{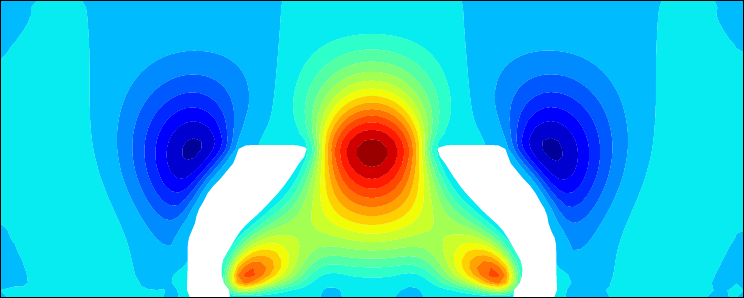}
        \subcaption{$z$ component at $t=2/3N_t$}
    \end{minipage}
    \begin{minipage}[t]{0.06\columnwidth}
        \centering
        \includegraphics[width=\columnwidth]{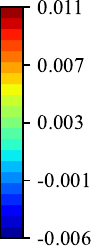}
    \end{minipage}
    \caption{Velocity distribution of 3D rotor. The optimized shape absorbs fluid from the side and pushes it out from the top. Here, $N_t = 3000\Delta t$.}
    \label{fig:example3_velocity}
\end{figure}
\begin{figure}
    \centering
    \includegraphics[width=0.5\columnwidth]{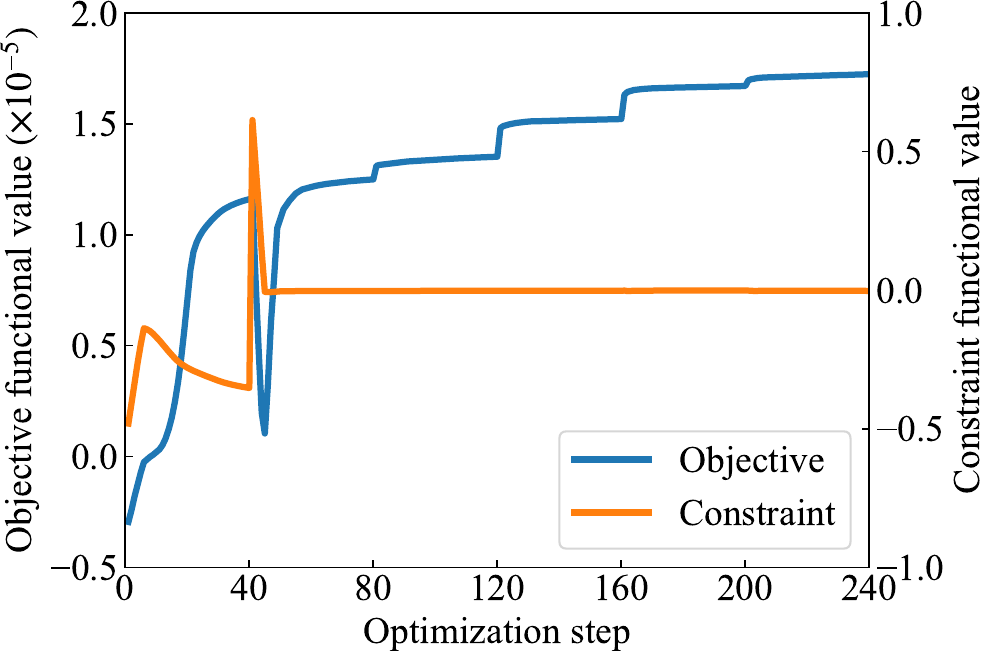}
    \caption{History plot of 3D rotor. The objective and constraint values change significantly near the $40$th step due to updating the prescribed volume limit.}
    \label{fig:example3_convergence}
\end{figure}
\begin{figure}
    \centering
    \includegraphics[width=0.5\columnwidth]{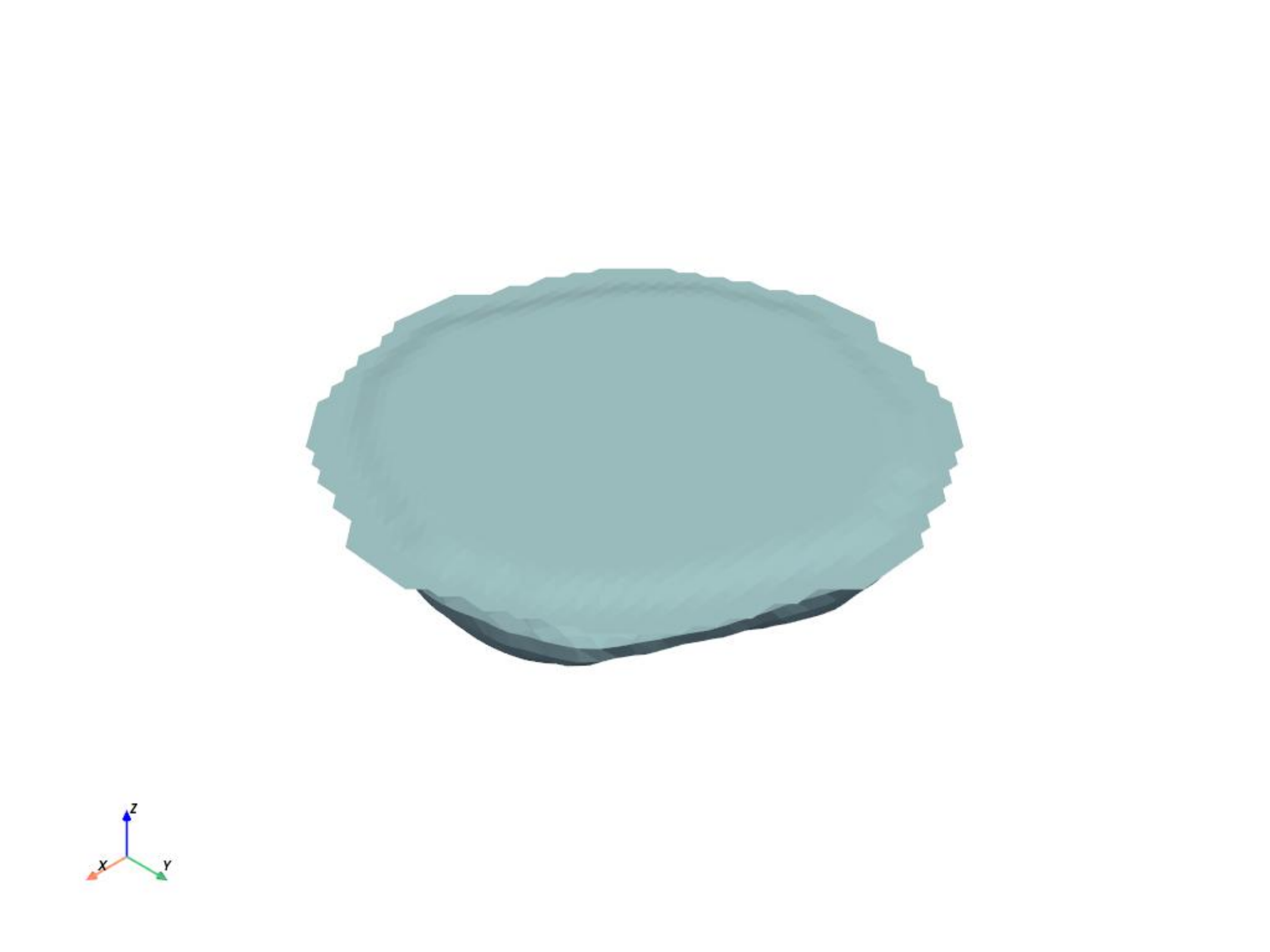}
    \caption{Poor local optimum of 3D rotor, in which there is no hole at the top.}
    \label{fig:example3_local_optimum}
\end{figure}
The optimized shape is shown in Fig.~\ref{fig:example3_optimized_shape}, and the distributions of the $x$-component and $z$-component of the fluid velocity on the plane $y=50\Delta x$ at each time are illustrated in Fig.~\ref{fig:example3_velocity}.
The optimized shape absorbs fluid from the side and push it out on the top.
It is noteworthy that, at $t=1/3N_t$, a temporary flow reversal occurs at the side of the rotor, as shown in Fig.~\ref{fig:example3_flow_x_t1000}.
Specifically, rather than absorbing fluid, the rotor briefly expels it.
This phenomenon arises because, at that moment, the rotor edge reaches the observation point, causing the fluid to separate.
This underscores the importance of treating the problem as an unsteady flow to accurately capture such transient behaviors.
The objective and constraint functional values over the course of the optimization are shown in Fig.~\ref{fig:example3_convergence}.
The objective functional exhibits a generally monotonic increase, except for a sudden drop observed shortly after the $40$th iteration.
This significant decline corresponds to the step at which the maximum allowable volume $V_{\text{max}}$ is updated from $100\%$ to $40\%$, as previously described.
At that point, the volume constraint is temporarily violated, but the shape is subsequently modified over a few iterations to satisfy the new constraint.

Initially, we also tested applying the volume constraint from the first iteration.
However, as shown in Fig.~\ref{fig:example3_convergence}, this constraint hindered the optimizer from exploring higher volume design, such as the peak observed between the 1st and 6th iterations in the unconstrained case.
As a result, the final shape lacked the top hole seen in Fig.~\ref{fig:example3_optimized_shape}, leading instead to the suboptimal configuration depicted in Fig.~\ref{fig:example3_local_optimum}, which differs significantly in topology.
The objective functional value for this alternative shape was markedly lower, indicating that it is a poor local optimum.
It is worth noting that the initial design field distribution is uniform, with $\gamma=0.5$, meaning the initial geometry contains no hole.
However, the optimized shape includes a hole to absorb and push out fluid, indicating that a topological change has occurred and plays a critical role in the design performance.
Therefore, topology optimization proves to be effective in this numerical example.  

Besides, the discontinuities in the objective functional value observed at the $80$th, $120$th, $160$th and $200$th iterations are due to the continuation process, where the steepness parameter $\beta$ of the projection filter is updated--consistent with the other numerical example. 

\section{Conclusion}
\label{sec4}

In this study, we propose the topology optimization method for moving objects in fluid.
The key idea of the proposed method is to separate grids for the analysis domain and the design domain.
The design variables and the design sensitivities are defined on the design grid, while the state and adjoint fields are defined on the analysis grid.
At each time step, the design grid first undergoes rigid body motion and is then mapped onto the analysis grid.
In the forward analysis, the Brinkman model coefficients and solid velocities are transferred from the design grid to the analysis grid, whereas in the reverse analysis, the state and adjoint fields are transferred in the reverse direction.
The accuracy of the design sensitivity computed using the proposed method was verified thorough comparison with a finite difference approximation.
The results confirmed that the proposed sensitivity analysis method achieves practically sufficient accuracy.
The proposed method was applied to three numerical example: a two-dimensional rotating case, a two-dimensional translating case, and a three-dimensional rotating case.
In each case, the optimized shape was physically interpretable and validated based on the flow characteristics.
In the future, this method will be extended to accommodate passive motions, where objects are moved by fluid flow rather than being actively driven.
Additionally, the framework holds promise for the co-design of moving and stationary components, such as turbines and their casings.
 
\section{Replication of results}
The necessary information for replication of the results is presented in this paper.
The interested reader may contact the corresponding author for further implementation details.

\section*{Acknowledgements}
This work was supported by JSPS KAKENHI (GrantNo. 23K26018).

\section{Verification of the representation of moving objects}

\begin{figure}[t]
    \centering
    \begin{minipage}[t]{0.3\columnwidth}
        \centering
        \includegraphics[width=0.8\columnwidth]{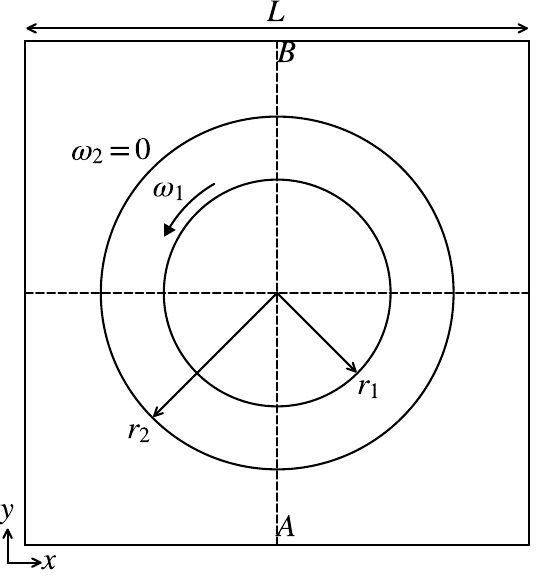}
        \subcaption{Analysis setting}
        \label{fig:verification_design_setting}
    \end{minipage}
    \begin{minipage}[t]{0.3\columnwidth}
        \centering
        \includegraphics[width=\columnwidth]{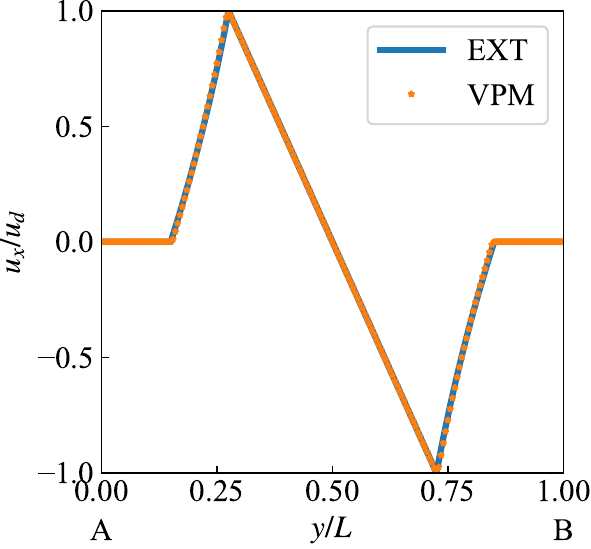}
        \subcaption{Comparison result}
        \label{fig:verification_comparison_result}
    \end{minipage}
    \begin{minipage}[t]{0.3\columnwidth}
        \centering
        \includegraphics[width=0.8\columnwidth]{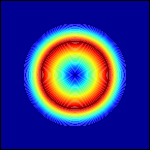}
        \subcaption{Velocity norm distribution}
        \label{fig:verification_velocity_distribution}
    \end{minipage}
    \begin{minipage}[t]{0.065\columnwidth}
        \centering
        \includegraphics[width=\columnwidth]{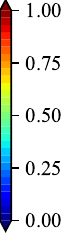}
    \end{minipage}
    \caption{Analysis setting and results for verification of the representation of the moving object. (a) The inner circle rotates with angler velocity $\omega_1$, while the outer circle does not rotate. (b) ``VPM'' means the values by the proposed method, while ``EXT'' means theoretical value. (c) The velocity norm distribution closely resembles that reported in the literature.}
\end{figure}

We verify the representation of moving objects using the proposed method.
The design setting is shown in Fig.~\ref{fig:verification_design_setting}, where $L=200\Delta x$, $r_1=45\Delta x$, $r_2=70\Delta x$ and $\omega_1=u_d/r_1$.
The results are compared with the theoretical value as follows:
\begin{equation}
    u_{\text{theoretical}}=\frac{-u_dr_1}{r_2^2-r_1^2}r+\frac{u_dr_1r_2^2}{r_2^2-r_1^2}\frac{1}{r},
\end{equation}
where $r$ is the distance from the center of rotation.
Fig.~\ref{fig:verification_comparison_result} shows a comparison along the dashed line from point A to point B in Fig.~\ref{fig:verification_design_setting}, demonstrating agreement between the numerical and theoretical results.
Furthermore, Fig.\ref{fig:verification_velocity_distribution} presents the flow velocity distribution, which closely resembles the distribution reported in the literature\cite{BENAMOUR2020101050}. 
\renewcommand{\thefigure}{\arabic{figure}} \setcounter{figure}{0}

\section*{Conflict of interest}
The authors declare that they have no conflict of interest.

\end{document}